\documentclass[seceqn,secthm]{elsart} 
 
\usepackage[applemac]{inputenc}  
\usepackage[english]{babel}
\usepackage{amsfonts,amsmath}
\usepackage{latexsym}
\usepackage{amssymb}
\journal{Stochastic Processes and their Applications}

\newcommand{\rz}[1]{\mathbb{#1}}
\def\R{\rz R}  
\def\N{\rz N}  
\def\shc{{\cal C}}
\def\shf{{\cal F}}
\def\shn{{\cal N}}

\def\indi{\mbox{1\hspace{-.25em}I}}

\begin{document}

\begin{frontmatter}
\title{Approximation via regularization of the local time of semimartingales and Brownian motion}
 \author[aut1]{B\'erard Bergery Blandine},  \author[aut1]{Vallois Pierre}
\address[aut1]{{\small{Universit\'e Henri Poincar\'e, Institut de
Math\'ematiques Elie Cartan, B.P. 239, F-54506 Vand\oe uvre-l\`es-Nancy Cedex, France}}}

\begin{abstract}
Through a regularization procedure, few approximation schemes of the local time of a large class of one dimensional processes are given.  We mainly consider the local time of continuous semimartingales and reversible diffusions, and the convergence holds in ucp sense. In the case of standard Brownian motion, we have been able to determine a rate of convergence in $L^2$, and a.s. convergence of some of our schemes.
\end{abstract}

\begin{keyword}
 local time \sep stochastic integration by regularization \sep quadratic variation \sep rate of convergence \sep stochastic Fubini's theorem\\
\textit{2000 MSC:} 60G44 \sep 60H05 \sep 60H99 \sep 60J55 \sep 60J60 \sep 60J65
\end{keyword}
\end{frontmatter}

\section{Introduction}\label{artsec1}
Let $(X_t)_{t\geqslant 0}$ be a continuous process which is defined on a  complete probability space  $(\Omega,   \shf,   \shf_{t},   P)$. It is supposed that  $(\shf_{t})$ verifies the usual hypotheses.

\textbf{1.} Let $X_t= M_t + V_t$ be a continuous $(\shf_{t})$-semimartingale,   where $M$ is a local martingale and $V$ is an adapted process with finite variation. In the usual stochastic calculus, two fundamental processes are associated with $X$: its quadratic variation and its local time. Indeed, It\^o's formula related to functions $f\in\shc^2$ is the following:
$$
f(X_t)= f(X_0) + \int_0^t f'(X_s) dX_s + \frac{1}{2}\int_0^t f''(X_s)d <X>_s,  
$$
where $<X>=<M>$ is the quadratic variation of either $X$ or $M$. 

  The random measure $ g \to \int_0^t g(X_s)d <X>_s$ is absolutly continuous with respect to the Lebesgue measure:  there exists  a measurable family of random variables $(L_t^x(X),   x\in \R,   t \geqslant 0)$, called the local time process related to $X$, such that
 for any non-negative Borel functions $g$:
\begin{equation}
\label{article02}
\int_0^t g (X_s)d<X>_s= \int_{\R} g(x) L_t^x(X) dx. 
\end{equation}
Moreover, $t\to L_t^x(X) $ is  continuous  and non-decreasing, for any $x\in \R$.

  Besides, an extension of It\^o's formula  in the case of convex functions $f$ may be obtained thanks to local time processes. Namely,
$$
f(X_t)=f(X_0) + \int_0^t f_{-}'(X_s)dX_s+ \frac{1}{2} \int_{\R}L_t^x f''(dx),  
$$
with $ f_{-}'$ the left derivative of $f$ and $f''$ the second derivative of $f$ in the distribution sense.

   \textbf{2.} The density occupation formula (\ref{article02}) gives a relation between the quadratic variation of a semimartingale and its local time process. We would like to show that the existence of the quadratic variation implies the weak existence of the local time process at a fixed level. For our purpose, it is convenient to use the definition of the quadratic variation which has been given in \cite{6}, in the setting of stochastic integration by regularization  \cite{8}: for any continuous process $X$,  its quadratic variation $ [X]_t$ is the process:
\begin{equation}
\label{article03b}
 [X]_t=\lim_{\epsilon \to 0} (ucp) \frac{1}{\epsilon}\int_0^t \left( X_{s+\epsilon}-X_{s}  \right)^2 ds,  
 \end{equation}
 provided that this limit exists in the (ucp) sense.  We  denote by (ucp) the convergence in probability,   uniformly on the compact sets (c.f. Section II.4 of \cite{16}). In the case of a continuous semimartingale $X$,   the quadratic variation defined by (\ref{article03b})  coincides with the usual quadratic variation $<X>$.  

  Let $\epsilon>0$  and $X$ be a continuous process. Let us introduce a family \\
  $(J_{\epsilon}(t,  y), y\in\R, t \geqslant 0)$ of processes, which will play a central role in our study:
\begin{equation}
\label{article05}
J_{\epsilon}(t,  y)=  \frac{1}{\epsilon}\int_0^t
 \left( \indi_{\{ y<X_{s+\epsilon}\}} - \indi_{\{ y<X_{s}\}} \right)  \left( X_{s+\epsilon}-X_{s}  \right)ds.
\end{equation}
It is actually possible (see \textbf{3}. of Section 2 for details) to prove that, if $X$ is a continuous process such that $[X]$ exists, then the measures $(J_{\epsilon}(t,  y)dy)$ on $\R$ weakly converge as $\epsilon \to 0$:
\begin{equation}
\label{article05c}
 \lim_{\epsilon \to 0} (ucp)\int_{\R}   f(y) J_{\epsilon}(t,  y) dy =\int_0^t f(X_s) d[X]_s,
\end{equation}
for any continuous function $f$ with compact support.

  As a consequence, if we suppose that $X$ is a semimartingale, the occupation times formula implies:
$$ \lim_{\epsilon \to 0} (ucp)\int_{\R} f(y) J_{\epsilon}(t,  y) dy  = \int_{\R} f(y) L_t^y(X) dy.$$
Thus,  in a certain sense, the measures $(J_{\epsilon}(t,  y)dy)$ converge to $(L_t^y(X) dy)$ as $\epsilon \to 0$. As a result, it seems natural to study  the convergence of $ J_{\epsilon}(t,  y) $ to $L_t^y(X) $ when $\epsilon \to 0$ and $y$ is a fixed real number. 

  Let us remark that, if $J_{\epsilon}(t,  y)$ converges in (ucp) sense, then its limit is equal to the covariation $[ X, \indi_{\{ y<X\}}]_t$ (c.f.  (\ref{article04a}) for the definition of the covariation).

   \textbf{3.} Our first approximation result concerns $(J_{\epsilon}(t,  y))$ when $X$ belongs to a class of diffusions stable under time reversal. This kind of diffusions has been studied in \cite{25} and \cite{26}. We consider the generalization made in Section 5 of \cite{24}. Let $X$ be a diffusion which satisfies
\begin{equation}
\label{article060}
X_t= X_0+ \int_0^t \sigma(s,X_s) dB_s + \int_0^t b(s,X_s)ds.
\end{equation}
It is moreover assumed that the following conditions hold:
\begin{equation}
\label{article060b}
\left.
\begin{array}{l}
  \forall t\in [0,T], X_t \textrm{ has a density } p(t,x) \textrm{ with respect to Lebesgue measure},       \\
      \sigma, b \textrm{ are jointly continuous},\\
        \sigma^2(s, .)\in W^{2,1}_{loc}(\R), b(s,.)\in W^{1,1}_{loc}(\R), xp(s,.)\in W^{2,\infty}_{loc}(\R), \\
        p\sigma^2 (s,.)\in W^{1,1}_{loc}(\R) \textrm{ hold for almost every } s\in [0,T],\\
      \frac{\partial p\sigma^2}{\partial x}, \frac{\partial^2 xp}{\partial x^2} \in L^1 ( [0,t]\times \R), t\in ]0,T].
\end{array}
\right\}
\end{equation}
To a fixed $T>0$, we associate the process
\begin{equation}
\label{articlereverse}
\widetilde{X}_u=X_{T-u},   \quad u\in [0,T].
\end{equation}
According to Theorem 5.1 of \cite{24}, $(\widetilde{X}_u)_{u\in [0 ,T]}$ is a diffusion which verifies the following equation:
\begin{equation}
\label{article060c}
\widetilde{X}_u=  X_T+  \int_0^u \sigma(T-s,\widetilde{X}_s) d\beta_s + \int_0^u \tilde{b}(T-s,\widetilde{X}_s)ds,
\end{equation}
where $\beta$ is a Brownian motion on a possibly enlarged space and $\tilde{b}$ is explicitly known. 

\begin{thm}
\label{article06}  Let $X$ be a diffusion which satisfies  (\ref{article060})-(\ref{article060b}). Then
$$  \lim_{\epsilon \to 0}(ucp)  \: J_{\epsilon}(t,  x)=L_t^x(X),\quad \forall x\in\R .$$
 \end{thm}  
Our limit in Theorem \ref{article06} is valid when $x$ is fixed. We have not been able to prove that the convergence is uniform with respect to $x$ varying in a compact set.

   \textbf{4.}  For simplicity of notation,   we take $x=0$ and we note $J_{\epsilon}(t)$ instead of $J_{\epsilon}(t,  0)$. The proof of Theorem \ref{article06}  is  based on a decomposition of   $J_{\epsilon}(t)$ as a sum of two terms. It is actually possible to prove (see Section \ref{artsec2}) that each term has a limit. However, theses limits cannot be expressed only through $L_t^0(X)$. Modifying the factor $ \indi_{\{ 0<X_{s+\epsilon}\}} - \indi_{\{ 0<X_{s}\} }$ in $J_{\epsilon}(t)$ and developing the product gives (see point \textbf{4}. of Section \ref{artsec1b} for details):
\begin{equation}
\label{article07a}
J_{\epsilon}(t) = I^3_\epsilon(t) + I^4_\epsilon(t)  + R_\epsilon(t), 
\end{equation}
where
\begin{eqnarray}
 I^3_\epsilon(t) & =&
 \frac{1}{\epsilon}\int_0^{t} X_{(u+\epsilon)\wedge t}^+ \indi_{ \{X_u \leqslant 0\} } du
+
\frac{1}{\epsilon}\int_0^{t} X_{(u+\epsilon) \wedge t}^- \indi_{ \{X_u>0\} } du,   \label{article07}\\
 I^4_\epsilon(t)&=&
  \frac{1}{\epsilon}\int_0^{t} X_u^- \indi_{ \{X_{(u+\epsilon) \wedge  t} > 0\} } du+  \frac{1}{\epsilon}\int_0^{t} X_u^+ \indi_{ \{X_{(u+\epsilon) \wedge t} \leqslant 0\} } du  \label{article08}, 
\end{eqnarray}
and  $(R_\epsilon(t))_{t\geqslant 0}$ is a process which goes  to 0 a.s. as $\epsilon \to 0$, uniformly  on compact sets in time.  Note that in (\ref{article07}) and (\ref{article08}), we have systemically introduced $X_{(u+\epsilon) \wedge  t} $ instead of $X_{u+\epsilon } $, in order to guarantee that $ I^3_\epsilon(t), I^4_\epsilon(t)$ are adapted processes.  This will play an important role in the proof of our results, in particular to obtain the convergence in the (ucp) sense.

  The decomposition (\ref{article07a}) seems at first complicated. However, it is interesting since it may be proved (see Theorem \ref{article10} below) that, under suitable assumptions, $ I^3_\epsilon(t)$ and $I^4_\epsilon(t)$ converge to a fraction of $L_t^0(X)$, as $\epsilon \to 0$.  
\begin{thm}
\label{article10}
\begin{itemize}
  \item[i)] If $X$ is a continuous semimartingale,   then
$$ \lim_{\epsilon \to 0}\: \textrm{ (ucp) } I^3_{\epsilon}(t)= \frac{1}{2} L_t^0(X). $$
  \item[ii)] If $X$ is a diffusion which satisfies (\ref{article060}) and (\ref{article060b}),   then
$$ \lim_{\epsilon \to 0}\: \textrm{ (ucp) } I^4_{\epsilon}(t)= \frac{1}{2} L_t^0(X). $$
\end{itemize}
\end{thm}
\begin{rem}
\begin{enumerate}
  \item We would like to emphasize that the choice of strict or large inequalities in (\ref{article07}) and (\ref{article08}) is important. Indeed,   the Lebesgue measure of $\{u; X_u=0\}$ may not vanish. Note that, in Tanaka's formula (c.f (\ref{arttanaka}) below), $X_t^+$ is associated with $ \indi_{ \{X_t>0\} }$. This heuristically explains our choice.
    \item   When $X$ is a standard Brownian motion, we will prove (see Theorem \ref{article11} below) that the two terms of the sum in $I^3_\epsilon(t)$ and $I^4_\epsilon(t)$ converge separately to $\frac{1}{4}L_t^0(X)$.
  \item Applying Theorem \ref{article10} to the process $(X_t-x)_{t\geqslant 0}$ gives  
  $$\lim_{\epsilon \to 0}\: \textrm{ (ucp) } I^3_{\epsilon}(t,x)= \frac{1}{2} L_t^x(X), \quad \textrm{ where }$$
  $$ I^3_\epsilon(t,x) =
 \frac{1}{\epsilon}\int_0^{t} (X_{(u+\epsilon)\wedge t}-x)^+ \indi_{ \{X_u \leqslant x\} } du
+
\frac{1}{\epsilon}\int_0^{t} (X_{(u+\epsilon) \wedge t}-x)^- \indi_{ \{X_u>x\} } du.$$
 There exists obviously a similar result where  $ I^3_{\epsilon}(t,x)$ is replaced by $ I^4_{\epsilon}(t,x)$.
\end{enumerate}
 \end{rem}

   \textbf{5.}  In the Brownian case, we will give complements to Theorem \ref{article06} and Theorem \ref{article10}. The first one concerns Theorem \ref{article10}. Obviously, $ I^3_\epsilon(t)$ and $ I^4_\epsilon(t)$  can be decomposed as:
\begin{equation}
\label{article0809}
I^3_\epsilon(t) =  I^{3,1}_\epsilon(t) +  I^{3,2}_\epsilon(t) +  r^3_\epsilon(t), \quad
I^4_\epsilon(t) =  I^{4,1}_\epsilon(t) +  I^{4,2}_\epsilon(t) + r^4_\epsilon(t),
\end{equation}
where
\begin{eqnarray}
 I^{3,1}_\epsilon(t)&=& \frac{1}{\epsilon}\int_0^{t} X_{(u+\epsilon) \wedge t}^- \indi_{ \{X_u>0\} }
 du, 
 \,
  I^{3,2}_\epsilon(t)   = \frac{1}{\epsilon}\int_0^{t} X_{(u+\epsilon)\wedge t}^+ \indi_{ \{X_u < 0\} }  du.\label{article09c}\\
  r^3_\epsilon(t) &=&  \frac{1}{\epsilon}\int_0^{t} X_{(u+\epsilon)\wedge t}^+ \indi_{ \{X_u = 0\} } du,\label{article09brestec}\\
 I^{4,1}_\epsilon(t)&=& \frac{1}{\epsilon}\int_0^{t} X_u^- \indi_{ \{X_{(u+\epsilon) \wedge  t} > 0\} } du, %
 \,
  I^{4,2}_\epsilon(t)   = \frac{1}{\epsilon}\int_0^{t} X_u^+ \indi_{ \{X_{(u+\epsilon) \wedge t} < 0\} } du,\label{article09b}\\
  r^4_\epsilon(t) &=&  \frac{1}{\epsilon}\int_0^{t} X_u^+ \indi_{ \{X_{(u+\epsilon) \wedge t}= 0\} } du,\label{article09breste}
  \end{eqnarray}
and $r^3_\epsilon(t), r^4_\epsilon(t)  $ converge a.s. to 0 as $\epsilon \to 0$ (c.f point \textbf{1}. of Section 5).
\begin{thm}
\label{article11}
Let $X$ be a standard Brownian motion. Then,
\begin{enumerate}
\item $$ \lim_{\epsilon \to 0}\: \textrm{ (ucp) } I^{3,1}_{\epsilon}(t)=  \lim_{\epsilon \to 0}\: \textrm{ (ucp) } I^{3,2}_{\epsilon}(t)= \frac{1}{4} L_t^0(X), $$
\item $$ \lim_{\epsilon \to 0}\: \textrm{ (ucp) } I^{4,1}_{\epsilon}(t)=  \lim_{\epsilon \to 0}\: \textrm{ (ucp) } I^{4,2}_{\epsilon}(t)= \frac{1}{4} L_t^0(X). $$
Moreover,   for all $T>0,  \delta \in ]0,  \frac{1}{2}[$,   there exists a constant $C$ such as
$$
 \left\| \sup_{t\in [0,  T]} \left| I^{4,i}_\epsilon(t) -   \frac{1}{4}L^0_t(X) \right| \right\|_{L^2(\Omega)} \leqslant  C \epsilon^{\frac{\delta}{2}},  \quad i= 1 \textrm{ or } 2.
$$
\end{enumerate}
\end{thm}
\begin{rem}
\label{article11rem}
The ucp convergence of $I^{4,1}_{\epsilon}(t), I^{4,2}_{\epsilon}(t)$ in Theorem \ref{article11} is still true (see \cite{0}) if $X_t=\int_0^t \sigma(s)dB_s$,   where $B$ is the standard Brownian motion and $\sigma: \R_+ \to \R$ a function which is H\"older continuous of order $\gamma > \frac{1}{4}$ and such as $| \sigma(s) | \geqslant  a > 0$.
\end{rem}
Then, we will give  below a complement related to Theorem \ref{article06}. In Theorem \ref{article06b}, we determine the rate of convergence of $J_{\epsilon}(t,  x) $ to $L_t^x(X)$ in $L^2(\Omega)$, as $\epsilon \to 0$. 
\begin{thm}
\label{article06b}     
Let $X$ be the standard Brownian motion.  For all $T>0,   x \in \R,   \delta \in ]0,  \frac{1}{2}[$,   there exists a constant $C$ such as:
$$ \forall \epsilon \in ]0,  1],   \quad \left\| \sup_{t\in [0,  T]} \left| J_{\epsilon}(t) -  L_t^0(X) \right|  \right\|_{L^2(\Omega)}  \leqslant C \epsilon^{\frac{\delta}{2}}. 
$$
\end{thm} 
In the setting of stochastic integration by regularization (see for instance \cite{5},   \cite{3},   \cite{4},    \cite{6} and \cite{8}), the ucp convergence is mainly used. So, Theorem \ref{article06}, \ref{article10}, \ref{article11} are of this type. It seems interesting to investigate the almost sure convergence. There exists few results using this type of convergence in \cite{15}. Our version of Theorem \ref{article06} formulated in terms of a.s. convergence is the following.
\begin{prop}
\label{article06bb}     
Let $X$ be the standard Brownian motion and $ (\epsilon_n)_{n\in \N}$ be a decreasing sequence of non-negative real numbers  which  satisfies  $\sum_{i=1}^{\infty} \sqrt{\epsilon_i}<\infty$.   Then, for any $x\in \R$, almost surely, 
\begin{eqnarray*}
\lim_{n\to \infty} \, \sup_{t\in [0,  T]} \left| J_{\epsilon_n}(t,  x) -  L_t^x (X)\right| & = & 0, \\
 \lim_{n\to \infty} \, \sup_{t\in [0,  T]} \left| I^{4,i}_{\epsilon_n}(t,x) - \frac{1}{4} L_t^x (X)\right| & = & 0, \quad {i= 1, 2}.
\end{eqnarray*}
\end{prop}

   \textbf{6.}  Let $(X_t)_{t\geqslant 0}$ be a continuous process  with values in $\R$. Let us briefly enumerate processes which admit local time processes $( L_t^x(X), x\in \R, t\geqslant 0)$. The most popular ones are semimartingales (see for instance Section VI of \cite{1}). Moreover, there exists a version of $( L_t^x(X), x\in \R, t\geqslant 0)$ such that $x\to L_t^x$ is right-continuous for any $t$. When $X$ is a Markov process, the local time process $(L_t^x)_{t\geqslant 0}$ at level $x$ is defined as a specific additive functional (see \cite{20}), and  through excursions in Chapter IV of \cite{21}. A construction  of the local time process related to one dimensional diffusions has been given in Chapter VI of \cite{13}. A large class of Lévy processes admits local time process (see Chapter V of  \cite{21}). Beyond semimartingales and Markov processes,  local time process associated with Gaussian processes may exist (c.f. \cite{19}) as occupation densities. A particular example of Gaussian process which has local time process is the fractional Brownian motion. Tanaka formula and It\^o-Tanaka formula have been given in \cite{23}.

   The definition of  $( L_t^x(X), x\in \R, t\geqslant 0)$ may not directly refer to the paths of  $(X_t)_{t\geqslant 0}$.  For instance, if $X$ is a continuous local martingale, $( L_t^x(X), t\geqslant 0)$ may be defined as the unique adapted process vanishing at 0 such that $| X_t-x| -L_t^x(X)$ is a local martingale. Therefore, it may be interesting to perform approximation schemes of $( L_t^x(X), x\in \R, t\geqslant 0)$ which involve more directly the paths of $X$. In the case of semimartingales, we deduce easily from the occupation times formula and right continuity of $t\to L_t^x(X)$ that $\frac{1}{\epsilon}\int_0^t \indi_{\{ x \leqslant X_s \leqslant x+\epsilon \}} d<X>_s$ tends to $L_t^x(X)$ as $\epsilon \to 0$.

   The Lévy excursion theory (c.f. \cite{14}) gives  other kind of approximation by the count of the downcrossings number or of the excursion number before a given time.  In the case of diffusion, the convergence of normalized sums to local time have been studied in \cite{11} and \cite{12}. Finaly, we may refer to  \cite{17} and  \cite{18}  for approximations of the local time for  Lévy's processes.

   \textbf{7.} Let us briefly detail the organization of the paper. Section \ref{artsec1b} contains few preliminary lemmas and complements related to some results stated in the Introduction. The proof of Theorem \ref{article06} is given in Section \ref{artsec2}.  Section \ref{artsec3} contains the proof of Theorem \ref{article10}. The Brownian case is studied in Section \ref{artsec4} (proof of Theorem \ref{article11}) and in  Section \ref{artsec5} (proof of Theorem \ref{article06b}). 

  Some results of this paper were announced without any proof in \cite{00}. Moreover, the setting in \cite{00} was the  Brownian one. 

  Let us adopt two conventions, which will be used in the sequel of the paper:
\begin{itemize}
\item  $[0,  T]$ will denote a given compact interval of time,
\item in the calculations, $C$ will denote a generic constant.
\end{itemize}

\section{Decomposition of $J_{\epsilon}(t)$ and preliminary lemmas }\label{artsec1b}
\setcounter{equation}{0}

This section has five independent parts. In the two first ones, we recall some known results related to Fubini's stochastic theorem (c.f. Lemma \ref{articlefubini} below) and H\"older continuity properties (c.f. Lemma \ref{article118} below). Proofs of developments given in the Introduction may be found in points \textbf{3}. and \textbf{4}. We end this Section by a technical lemma.

   \textbf{1.} Let us start with a modification of Fubini's theorem. This result may be found in Section IV.5 of  \cite{1} and is  crucial in most of our proofs. It permits to express some Lebesgue integrals as stochastic integrals with respect to martingales. This allows to obtain (ucp) convergence via Doob's inequality (see for instance the proof of Proposition \ref{article064a}).
\begin{lem} 
\label{articlefubini}
Let $(H(u,  s), s \in \R,  u \geqslant 0)$ be a collection of predictable processes which are measurable with respect to $(u,s,\omega)$. If one of the following hypotheses holds:
\begin{description}
\item[i.]  $(X_s)_{s \geqslant 0}$ is the standard Brownian motion and  $ \int_0^t \int_0^t E( H(u,  s)^2 ) ds du < \infty$,  
\item[ii.] $(X_s)_{s \geqslant 0}$ is a continuous semimartingale and $(H(u,  s))_{s,  u \geqslant 0}$ is  uniformly bounded,  
\end{description}
then,   almost surely,   
$$
\int_0^{t} \left[ \int _0^{t} H(u,  s)dX_s \right] du = \int _0^{t}  \left[ \int_0^{t} H(u,  s)du\right] dX_s,\quad \forall t \geqslant 0 . 
$$
\end{lem}

   \textbf{2.} Let us recall  some H\"older continuity  properties of the Brownian motion and its local time process.     
\begin{lem}
\label{article118}
Let $\delta \in ]0,   \frac{1}{2}[, T>0$  and $X$  be the standard Brownian motion. 
\begin{itemize}
  \item[i)] Then,  there exists a positive random  constant  $C_{\delta}\in L^2(\Omega)$ such as a.s.
$$
 | X_y-X_{y'}| \leqslant  C_{\delta} |y-y'|^{\delta}, \quad \forall y,   y' \in [0,   T]  .  
$$
  \item[ii)] The Brownian local time is H\"older continuous in space:  there exists a positive random  constant $K_{\delta}\in L^2(\Omega)$ such as a.s.  
$$ | L_t^{a}(X)- L_t^{a'}(X) | \leqslant K_{\delta} | a-a'|^{\delta}, \quad \forall  t \in [0,  T],   a,  a' \in \R. $$ 
\end{itemize}
\end{lem}
In the further proofs, as soon as $\delta \in ]0,   \frac{1}{2}[$ and $T$ are given,  the constants $K_{\delta}, C_{\delta}$ may be considered as fixed.

  \textbf{3.}  Let  $Y$ and $Z$ be continuous processes.  In \cite{4}, the covariation of $Y$ and $Z$ has been defined as:
\begin{equation}
\label{article04a}
[ Y,   Z]_t= \lim_{\epsilon \to 0} (ucp) \frac{1}{\epsilon}\int_0^t  \left( Y_{s+\epsilon}-Y_{s}  \right) \left( Z_{s+\epsilon} -Z_{s}  \right) ds,
\end{equation}
 if the limit exists.

  Let $X$ be a continuous process with finite quadratic variation $[X]=[X,X]$ (c.f.(\ref{article03b})). According to Proposition 2.1 of \cite{4}:
\begin{equation}
\label{article04}
 [f(X),  g(X)]_t=\int_0^t f'(X_s)g'(X_s)d[X]_s, \quad \forall f,  g\in \shc^1.
\end{equation}
The aim of this section is to explain how (\ref{article04}) combined with (\ref{article04a}) leads to (\ref{article05c}). Let $f$ be a continuous function  with compact support  and let $F$ be its primitive which vanishes at 0.  Taking $Y_t= F(X_t)$ and $Z_t= X_t$ in  (\ref{article04a}) and using (\ref{article04}) comes to:
\begin{eqnarray}
 \lim_{\epsilon \to 0} \, (ucp) \,  \frac{1}{\epsilon}\int_0^t \left( F(X_{s+\epsilon})-F(X_{s})  \right) \left( X_{s+\epsilon}-X_{s}  \right) ds &= &[F(X),  X]_t , \nonumber\\
 &=&\int_0^t f(X_s) d[X]_s. \label{article04c}
\end{eqnarray}
We now express the integral on the left hand-side of (\ref{article04c}) through $f$ instead of $F$. Since $F(x)=  \int_{\R} \left( \indi_{\{ y<x\}} - \indi_{\{ y<0\}} \right)f(y) dy$, applying Fubini's theorem gives:
$$  \frac{1}{\epsilon}\int_0^t \left( F(X_{s+\epsilon})-F(X_{s})  \right) \left( X_{s+\epsilon}-X_{s}  \right) ds =\int_{\mathbb{R}}J_{\epsilon}(t,y) f(y) dy.$$
It is now clear that (\ref{article05c}) is a direct consequence of (\ref{article04c}).

   \textbf{4.} We would like to prove (\ref{article07a}).  First, let us introduce $R_\epsilon(t)$:
\begin{eqnarray*}
 R_\epsilon(t)& = &  \frac{1}{\epsilon}\int_{(t-\epsilon)^+}^t
 \left( \indi_{\{ 0<X_{u+\epsilon}\}} - \indi_{\{ 0<X_{u}\}} \right)  \left( X_{u+\epsilon}-X_{u}  \right)du  \\
& & -   \frac{1}{\epsilon}\int_{(t-\epsilon)^+}^t
 \left( \indi_{\{ 0<X_{t}\}} - \indi_{\{ 0<X_{u}\}} \right)  \left( X_{t}-X_{u}  \right)du,  \quad t \geqslant 0.
\end{eqnarray*}
Then, 
\begin{equation}
\label{article07c}
J_{\epsilon}(t) -   R_\epsilon(t) =  \frac{1}{\epsilon}\int_0^{t}
( \indi_{\{ 0<X_{(u+\epsilon)\wedge t }\}} - \indi_{\{ 0<X_{u}\}} )  \left( X_{(u+\epsilon)\wedge t }-X_{u}  \right)du. 
\end{equation}
Using
\begin{equation}
\label{article07cbis}
 \indi_{\{X_{(u+\epsilon)\wedge t } >0\}} -\indi_{\{X_{u} >0\} }= \indi_{\{X_{(u+\epsilon)\wedge t  } >0,   X_u \leqslant 0\} }-\indi_{\{X_{(u+\epsilon)\wedge t } \leqslant 0,   X_u > 0\} },
\end{equation}
and developing the product of the integrand in (\ref{article07c}) lead to
\begin{eqnarray*}
J_{\epsilon}(t) -   R_\epsilon(t) &  =  &
\frac{1}{\epsilon}\int_0^{t}  \left[
X_{(u+\epsilon)\wedge t }\indi_{\{X_{(u+\epsilon)\wedge t } >0,   X_u \leqslant 0\} }
-X_{(u+\epsilon)\wedge t }\indi_{\{X_{(u+\epsilon)\wedge t } \leqslant 0,   X_u > 0\} }\right.  \\
& & \left.  \qquad
-X_u\indi_{\{X_{(u+\epsilon)\wedge t } >0,   X_u \leqslant 0\} }
+X_u\indi_{\{X_{(u+\epsilon)\wedge t } \leqslant 0,   X_u > 0\} } \right] du. 
\end{eqnarray*}
Since $X_{. }\indi_{\{X_{. } >0\} }=X_{. }^+$ and $-X_{. }\indi_{\{X_{. } \leqslant 0\} }=X_{. }^-$,  (\ref{article07a}) follows.

   \textbf{5.}  In (\ref{article07a}), $R_\epsilon(t)$ may be viewed as a remainder term. The lemma below  ensures the convergence to 0, in the a.s. sense, of this kind of terms. A proof of Lemma \ref{articlerestezero} may be found in \cite{0}.
\begin{lem}\label{articlerestezero}
Let $X$ be a continuous process and let $ a,  b,  c,  d : [0,  1]\times [0,  T+1]\rightarrow \R,  $ be Borel functions such that 
$$
 0 \leqslant b(\epsilon,  t ) -a (\epsilon,  t ) \leqslant \epsilon,   \quad |c(\epsilon,  s ) -d (\epsilon,  s )|  \leqslant \epsilon, \quad t,s \in[0,  T+1], \epsilon \in [0,1]. 
$$
Let us consider
 $$\widehat{R}_{\epsilon}(t)= \frac{1}{\epsilon}\int_{a(\epsilon,  t)}^{b(\epsilon,  t)} (X_{c(\epsilon,  s )}-X_{d(\epsilon,  s )}) \indi_{A_s} ds, \quad 0 \leqslant t \leqslant T, $$
where $(A_s)_{ s \in [0,  T+1]}$ is a collection of measurable events. \\
Then $\widehat{R}_{\epsilon}(t)$ tends a.s.  to 0 as $\epsilon \to 0$, uniformly on $[0,  T]$.  Furthermore,   if $X$ is the standard Brownian motion and $\delta\in ]0, \frac{1}{2}[$,  there exists a positive random  constant  $C_{\delta}\in L^2(\Omega)$ such that, almost surely, 
$$ \sup_{t\in [0,  T]} | \widehat{R}_{\epsilon}(t) | \leqslant  C_\delta \epsilon^\delta.$$
\end{lem}
%

\section{Proof of Theorem \ref{article06} and associated results}\label{artsec2}
\setcounter{equation}{0}

   \textbf{1.} 
Let $X$ be a continuous process. We have the following decomposition of $J_{\epsilon}(t)$:
\begin{equation}\label{article061}
J_{\epsilon}(t)=  -I^1_{\epsilon}(t)+I^2_{\epsilon}(t),  
\end{equation}
with
\begin{eqnarray}
I^1_{\epsilon}(t) & = & \int_0^t \frac{X_{s+\epsilon}-X_{s}}{\epsilon} \indi_{\{ 0<X_{s}\}} ds,  \label{article062} \\
I^2_{\epsilon}(t) & = & \int_0^t \frac{X_{s+\epsilon}-X_{s}}{\epsilon} \indi_{\{ 0<X_{s+\epsilon}\}} ds.  \label{article063}
\end{eqnarray}
The goal of this section is the study of the convergence of $I^1_{\epsilon}(t) $ and $I^2_{\epsilon}(t) $ as $\epsilon \to 0$.

  The first result concerns $I^1_{\epsilon}(t)$ (see point \textbf{2}. for the proof).
\begin{prop}
\label{article064a}
Let $X=M+V$ be a continuous semimartingale, with $M$ a continuous local martingale and $V$ an adapted continuous process with bounded variation. It is moreover supposed that both $dV$ and $d<M>$ are absolutely continuous with respect to the Lebesgue measure. Then
$$
\lim_{\epsilon \to 0}(ucp)  \: I^1_{\epsilon}(t)=  \int_0^t \indi_{\{ 0<X_{s}\}} dX_s.
$$
\end{prop}
\begin{rem}
If $X$ is a diffusion which satisfies (\ref{article060}), then Proposition \ref{article064a} applies and $I^1_{\epsilon}(t)$ converges to $ \int_0^t \indi_{\{ 0<X_{s}\}} dX_s$ in (ucp) sense when $\epsilon \to 0$.
\end{rem}

  Let $\widetilde{X}$ be the time reversal process, defined by (\ref{articlereverse}). The study of $I^2_{\epsilon}(t)$ may be reduced to the one of $\int_0^t \frac{\widetilde{X}_{s+\epsilon}-\widetilde{X}_{s}}{\epsilon} \indi_{\{ 0<\widetilde{X}_{s}\}} ds$ (c.f. point \textbf{3}. for details). 
Obviously, this term is of $I^1_{\epsilon}(t)$-type. Consequently, its convergence may be obtained by using Proposition \ref{article064a},  if we know that  $\widetilde{X}$ is a nice semimartingale. This property holds when $X$ is a diffusion which satisfies  (\ref{article060}) and  (\ref{article060b}). This justifies the interest of reversible diffusions.
\begin{prop}
 \label{article064b} 
If $X$ is a diffusion which satisfies  (\ref{article060}) and  (\ref{article060b}), then
$$
\lim_{\epsilon \to 0}(ucp)  \: I^2_{\epsilon}(t) =  X_t^+ - X_0^+  + \frac{1}{2}L_t^0(X). 
$$
\end{prop}
The proof of this result is postponed in point \textbf{3}. below.
\begin{rem}\begin{enumerate}
  \item When $X$ is a diffusion which satisfies  (\ref{article060}) and  (\ref{article060b}), then $I^1_{\epsilon}(t)$ and $I^2_{\epsilon}(t)$ converge as $\epsilon \to 0$.  Theorem \ref{article06} is a direct consequence of  Tanaka's formula: 
\begin{equation}
\label{arttanaka}
X_t^+=  X_0^+ +  \int_0^t \indi_{\{ 0< X_s  \}}dX_s + \frac{1}{2} L_t^0(X).
\end{equation}
  \item When $X$ is a standard Brownian motion, it can be proved directly that $J_{\epsilon}(t)$ converges to $L_t^0(X)$ (c.f. \cite{00}).
\end{enumerate}
\end{rem}

   \textbf{2.  Proof of Proposition \ref{article064a}.} Let $X$ be a continuous semimartingale with canonical decomposition $X=M+V$. The key of our proof is to write $I^1_{\epsilon}(t)$ as the sum of a semimartingale plus a remainder term. From a technical point of view, the semimartingale will be obtained by expressing $ I^1_{\epsilon}(t)$ through $X_{(s+\epsilon)\wedge t}$ instead of $X_{s+\epsilon}$. Namely,
\begin{equation}
\label{article068}
 I^1_{\epsilon}(t)-\int_0^t \indi_{\{ 0<X_{s}\}} dX_s = \widetilde{I}^1_{\epsilon}(t) + \widehat{I}^1_{\epsilon}(t) + \Delta_1(t,  \epsilon),  
\end{equation}
with
\begin{eqnarray}
\widetilde{I}^1_{\epsilon}(t) &=& \int_0^{t} \frac{1}{\epsilon}(M_{(s+\epsilon)\wedge t}-M_{s}) \indi_{\{ 0<X_{s} \}} ds -\int_0^t  \indi_{\{ 0<X_{s} \}}dM_s    \label{article069} \\
\widehat{I}^1_{\epsilon}(t) &=& \int_0^{t} \frac{1}{\epsilon}(V_{(s+\epsilon)\wedge t}-V_{s}) \indi_{\{ 0<X_{s} \}} ds -\int_0^t  \indi_{\{ 0<X_{s} \}}dV_s, \label{article069b}  \\
\Delta_1(t,  \epsilon)&= & \frac{1}{\epsilon} \int_{(t-\epsilon)^+}^t  \left( X_{s+\epsilon} -X_t \right) \indi_{\{ 0<X_{s} \}}ds.   \label{article0610}
\end{eqnarray}
By Lemma \ref{articlerestezero},   $\Delta_1(t,  \epsilon)$ tends a.s.  to 0,   uniformly on the compact sets. We will prove the convergence of  $\widehat{I}^1_{\epsilon}(t) $ (resp. $\widetilde{I}^1_{\epsilon}(t)$) in step a) (resp. b)) below.

  \textbf{a)} First, we will show that  $\widehat{I}^1_{\epsilon}(t) $  tends to 0 almost surely. The key of our approach is to write $\widehat{I}^1_{\epsilon}(t) $ as an integral with respect to $dV$, through Fubini's theorem:
\begin{eqnarray*}
\widehat{I}^1_{\epsilon}(t) &= &\int_0^{t}  \frac{1}{\epsilon} \left(  \int_s^{(s+\epsilon)\wedge t} \indi_{\{ 0<X_{s} \}}  dV_u\right)  ds -\int_0^t  \indi_{\{ 0<X_{s} \}}dV_s, \\
&=&\int_0^t    \left(  \frac{1}{\epsilon} \int_{(u-\epsilon)^+}^u  \indi_{\{ 0<X_{s}\}} ds - \indi_{\{ 0<X_{u}\}}\right) dV_u. 
\end{eqnarray*}
Thus 
$$
 \sup_{0\leqslant t \leqslant T} | \widehat{I}^1_{\epsilon}(t) |   \leqslant \int_0^T    \left|  \frac{1}{\epsilon} \int_{(u-\epsilon)^+}^u  \indi_{\{ 0<X_{s}\}} ds - \indi_{\{ 0<X_{u}\}}\right| d|V|_u.
$$
We observe that
\begin{equation}
\label{article0610b}
\lim_{\epsilon \to 0} \frac{1}{\epsilon} \int_{(u-\epsilon)^+}^{u} \indi_{\{ 0<X_{s} \}}ds -\indi_{\{ 0<X_{u} \}}=0, \quad (du)\textrm{ almost everywhere.} 
 \end{equation}
Since $\frac{1}{\epsilon} \int_{(u-\epsilon)^+}^{u} \indi_{\{ 0<X_{s} \}}ds$ is bounded by 1, and $d|V|_s$ absolutely continuous  with respect to the Lebesgue measure, then Lebesgue's convergence theorem implies that $\widehat{I}^1_{\epsilon}(t) $  tends to 0, uniformly for $t\in [0,T]$ , when $\epsilon\to 0$.

  \textbf{b)} Next, we will show that $\sup_{t\in [0,T]} \widetilde{I}^1_{\epsilon}(t) $ tends to 0 in $L^2(\Omega)$.  We will prove that $\widetilde{I}^1_{\epsilon}(t) $ is a stochastic integral.  Since $  M_{(s+\epsilon)\wedge t}-M_{s} = \int_s^{(s+\epsilon)\wedge t} dM_u  $ and  $\indi_{\{ 0<X_{s} \}} $ is adapted, we have:  
$$\widetilde{I}^1_{\epsilon}(t) = \int_0^{t}  \frac{1}{\epsilon} \left(  \int_s^{(s+\epsilon)\wedge t} \indi_{\{ 0<X_{s} \}}  dM_u\right)  ds -\int_0^t  \indi_{\{ 0<X_{s} \}}dM_s . $$
Stochastic Fubini's theorem (i.e. Lemma \ref{articlefubini}) may be applied:
$$\widetilde{I}^1_{\epsilon}(t)  = \int_0^t    \left(  \frac{1}{\epsilon} \int_{(u-\epsilon)^+}^u  \indi_{\{ 0<X_{s}\}} ds - \indi_{\{ 0<X_{u}\}}\right) dM_u. $$
As a result,   $\widetilde{I}^1_{\epsilon}(t)$ is a local martingale.

  Let us suppose that $<M>_T$ is bounded. It is clear that
\begin{eqnarray*}
<\widetilde{I}^1_{\epsilon}>_t  &=& \int_0^t \left( \frac{1}{\epsilon} \int_{(u-\epsilon)^+}^{u} \indi_{\{ 0<X_{s} \}}ds -\indi_{\{ 0<X_{u} \}} \right)^2 d<M>_u,\\
& \leqslant & \int_0^T 4 d<M>_u = 4 <M>_T
\end{eqnarray*}
Consequently, $(\widetilde{I}^1_{\epsilon})_{t\in [0,T]}$ is a martingale which is bounded in $L^2(\Omega)$. This allows to apply  Doob's inequality:
\begin{eqnarray}
 E\left(  \sup_{0\leqslant t \leqslant T}( \widetilde{I}^1_{\epsilon}(t) )^2 \right)& \leqslant&   4E\left(\left(\widetilde{I}^1_{\epsilon}(T) \right)^2\right),     \label{article0611} \\
& \leqslant & 4 E \left[\int_0^T  \left(\frac{1}{\epsilon} \int_{(u-\epsilon)^+}^{u}  \indi_{\{ 0<X_{s}\}}ds - \indi_{\{ 0<X_{u}\}}\right)^2 d<M>_u\right]. \nonumber
\end{eqnarray}
%
By using (\ref{article0610b}) and by repeating the arguments developed in item \textbf{a)} above, we may conclude that $E\left(  \sup_{0\leqslant t \leqslant T}( \widetilde{I}^1_{\epsilon}(t) )^2 \right) $ goes to 0 when $\epsilon\to 0$. This implies that $\lim_{\epsilon \to 0} \, (ucp)\,  \widetilde{I}^1_{\epsilon}(t)  = 0$.

  Since $<M>_T$ is not necessarly bounded, let us introduce the following sequence of stopping times:
 $$T_n = \inf\{t \geqslant 0, <M>_t \geqslant n\},$$
  with the convention $\inf \emptyset =+\infty$. Since $M^{T^n}$ is a local martingale so that \\
  $<M^{T^n}>_t= <M>_{t\wedge T_n} \leqslant n$, then, for any $n \geqslant 0$, $ \widetilde{I}^1_{\epsilon}(t\wedge T_n)$ goes to 0 as $\epsilon \to 0$, in the (ucp) sense. Using the fact that $T_n$ is a non decreasing sequence of stopping time converging to $+\infty$ as $n\to +\infty$, it follows that $\lim_{\epsilon \to 0} \textrm{ (ucp) }  \widetilde{I}^1_{\epsilon}(t) = 0$.\qed

   \textbf{3.  Proof of Proposition \ref{article064b}.} 
Let us now suppose that $X$ is a diffusion which satisfies (\ref{article060}) and (\ref{article060b}) . Since  $( \indi_{\{ 0<X_{s+\epsilon}\}})_{s\geqslant 0}$ is not a  $(\shf_s)$-adapted process,   we are not allowed to use  Fubini's theorem  as in point \textbf{2}. above.  We will use time reversal.  Roughly speaking, time reversal allows to reduce the study of $I^2_{\epsilon}(t)$ to the one of $I^1_{\epsilon}(t)$.  First,   we make the change of variable $u=T-s-\epsilon$ and we write $I^2_{\epsilon}(t)$ as a sum of two terms:
\begin{equation}
\label{article065} 
I^2_{\epsilon}(t)= \widetilde{I}^2_{\epsilon}(t) +\Delta_2(\epsilon,   t),  
\end{equation}
with  
\begin{eqnarray*}
 \widetilde{I}^2_{\epsilon}(t)& = &\indi_{\{ \epsilon \leqslant t \}}  \frac{1}{\epsilon}  \int_{T-t}^{T-\epsilon}(X_{T-u}-X_{T-u-\epsilon}) \indi_{\{ 0<X_{T-u}\}} du,   \\
\Delta_2(\epsilon,   t) & = & \frac{1}{\epsilon} \int^{T-(t\vee \epsilon)}_{T-t-\epsilon} (X_{T-u}-X_{T-u-\epsilon}) \indi_{\{ 0<X_{T-u}\}} du. 
\end{eqnarray*}
By Lemma \ref{articlerestezero},   $\Delta_2(t,  \epsilon)$ tends a.s.  to 0,   uniformly on the compact sets.  Thus,   the convergence of $ \widetilde{I}^2_{\epsilon}(t)$ implies the one of $I^2_{\epsilon}(t)$ to the same limit.  \\
Let us recall that $\widetilde{X}$ has been defined by (\ref{articlereverse}). By Theorem 5.1 of \cite{24}, $(\widetilde{X}_u)_{u\in [ 0 ,T ]}$ is a diffusion which satisfies (\ref{article060c}). Let us note that $<\widetilde{X}>_t = \int_0^t \sigma^2(T-s,\widetilde{X}_s)ds$. Then,
$$\widetilde{I}^2_{\epsilon} (t)= - \indi_{\{ \epsilon \leqslant t \}} \frac{1}{\epsilon}  \int_{T-t}^{T-\epsilon} ( \widetilde{X}_{u+\epsilon} -\widetilde{X}_u)\indi_{\{ 0< \widetilde{X}_u \}} du.$$
Proposition  \ref{article064a} yields to
$$ \lim_{\epsilon \to 0}(ucp)  \: \widetilde{I}^2_{\epsilon} (t)=- \int_{T-t}^{T}  \indi_{\{ 0 <\widetilde{X}_u \}} d\widetilde{X}_u. $$
We now express the limit of $\widetilde{I}^2_{\epsilon} (t)$ in terms of $L_t^0(X)$. By Tanaka's formula, we get:
$$\widetilde{X}_T^+ - \widetilde{X}_{T-t}^+  =  \int_{T-t}^T \indi_{\{ 0< \widetilde{X}_u \}} d\widetilde{X}_u  + \frac{1}{2} (L_T^{0}( \widetilde{X})- L_{T-t}^{0}( \widetilde{X})).$$
Since $\widetilde{X}_T= X_0$ and $\widetilde{X}_{T-t} = X_t$, we have
$$ \lim_{\epsilon \to 0}(ucp) \: \widetilde{I}^2_{\epsilon}(t)= X_t^+ - X_0^+ +  \frac{1}{2} (L_T^{0}( \widetilde{X})- L_{T-t}^{0}( \widetilde{X})) .
$$
In a last step, we express $(L_t^{0}( \widetilde{X})_{t\in [0,T]}$ via $(L_t^{0}( X))_{t\in [0,T]}$. Since $ x \to L_t^{x}( \widetilde{X})$ is right-continuous, we obtain:
\begin{eqnarray} 
L_T^0(\widetilde{X}) -  L_{T-t}^0 (\widetilde{X}) & = &
 \lim_{\alpha \to 0} \frac{1}{\alpha} \int_{T-t}^T \indi_{\{ 0<    \widetilde{X}_u <  \alpha \}} d<\widetilde{X}>_u, \nonumber\\
& = &  \lim_{\alpha \to 0} \frac{1}{\alpha}  \int_{T-t}^T \indi_{\{ 0<    X_{T-u} <  \alpha \}} \sigma^2(T-u, X_{T-u})du,  \nonumber\\
 &=& \lim_{\alpha \to 0} \frac{1}{\alpha} \int_0^t \indi_{\{ 0 < X_s <  \alpha \}}   \sigma^2(s, X_s)ds \quad ( s=T-u ), \nonumber  \\
 & = & \lim_{\alpha \to 0} \frac{1}{\alpha} \int_0^t \indi_{\{ 0 < X_s <  \alpha \}} d<X>_s = \frac{1}{2} L_t^0(X). \label{article065b}  
\end{eqnarray}
As a result,    $ \widetilde{I}^2_{\epsilon}(t)$ tends in the ucp sense to $X_t^+ -X_0^+   +  \frac{1}{2} L_t^0(X)$ when $\epsilon \to 0$.  
\qed

\section{Proof of Theorem \ref{article10}  }\label{artsec3}
\setcounter{equation}{0}

We have already observed in the proof of Theorem \ref{article06} (see Section 3) that time reversal property allows to reduce the convergence of $I^2_\epsilon(t)$ to the one of $I^1_\epsilon(t)$. The same property permits again to obtain the convergence of $I^4_\epsilon(t)$ via the one  of  $I^3_\epsilon(t)$. We begin with the study of $I^3_\epsilon(t)$ in point \textbf{1}, and we will deduce the convergence of  $I^4_\epsilon(t)$ in point \textbf{2}.

   \textbf{1. Proof of point \textit{i)}.} 
In this part,   $X$ will be a continuous semimartingale with canonical decomposition $X=M+V$.  Reasoning by stopping (see point \textbf{2.b} of the proof of Proposition \ref{article064a}) allows to suppose that $M$, $<M>$ and the total variation of $V$ are bounded processes.  Let us  recall that $ I^3_\epsilon(t)$ has been defined by (\ref{article07}). Our aim is to prove that $ I^3_\epsilon(t)$ goes to $\frac{1}{2}L_t^0$ as $\epsilon \to 0$.  Our approach is based mainly  on Tanaka's formula and Fubini's theorem.

   \textbf{1.a.} By using Tanaka's formula,  we get: 
$$X_{(u+\epsilon)\wedge t }^+=X_u^+ +\int_u^{(u+\epsilon)\wedge t }\indi_{\{X_s > 0\} } dX_s + \frac{1}{2}(L_{(u+\epsilon)\wedge t }^0(X)-L_u^0(X)). $$
Since $X_u^+ \indi_{ \{X_{u} \leqslant 0\} } du=0$, integrating the previous relation over $[0,T]$ gives
\begin{eqnarray*}
\frac{1}{\epsilon}\int_0^{t}  X_{(u+\epsilon)\wedge t }^+ \indi_{ \{X_u \leqslant 0\} } du &= &
\frac{1}{\epsilon}\int_0^{t} \left[  \indi_{ \{X_u \leqslant 0\} }\int_u^{(u+\epsilon)\wedge t }\indi_{\{X_s > 0\} } dX_s \right.  \\
&&\left.  + \frac{1}{2}(L_{(u+\epsilon)\wedge t }^0(X)-L_u^0(X)) \indi_{\{ X_u \leqslant 0\} } \right] du. 
\end{eqnarray*}
The process   $(\indi_{ \{X_u \leqslant 0\} } )_{u\geqslant 0}$ is adapted,   thus
\begin{eqnarray}
\frac{1}{\epsilon}\int_0^{t} X_{(u+\epsilon)\wedge t }^+ \indi_{ \{X_u \leqslant 0\} } du & =&
\frac{1}{\epsilon}\int_0^{t} \left( \int_u^{(u+\epsilon)\wedge t } \indi_{\{X_s > 0\} } \indi_{\{ X_u \leqslant 0\} }dX_s \right) du \nonumber \\
&&+
\frac{1}{2\epsilon}\int_0^{t} (L_{(u+\epsilon)\wedge t }^0(X)-L_u^0(X)) \indi_{\{ X_u \leqslant 0\} } du. \label{article10eq2}
\end{eqnarray}
The same method applies  to $X^-$ instead of $X^+$. Combining the two expressions comes to:
\begin{eqnarray*}
I^3_\epsilon(t) & = &  \frac{1}{\epsilon}\int_0^{t} \left( \int_u^{(u+\epsilon)\wedge t }\indi_{\{X_s > 0,   X_u \leqslant 0\} }dX_s\right) du \\
& & -  \frac{1}{\epsilon}\int_0^{t} \left(  \int_u^{(u+\epsilon)\wedge t }\indi_{\{X_s \leqslant 0,   X_u > 0\} } dX_s \right) du\\
& & +  \frac{1}{2\epsilon}\int_0^{t} (L_{(u+\epsilon)\wedge t }^0(X)-L_u^0(X))du. 
\end{eqnarray*}

   \textbf{1.b.  Study of $ \frac{1}{2\epsilon}\int_0^{t} (L_{(u+\epsilon)\wedge t }^0(X)-L_u^0(X))du$. }  \\
First,   we write $(L_{(u+\epsilon)\wedge t }^0(X)-L_u^0(X))$ as a Stieltjes integral with respect to $dL_.^0(X)$:
$$
\frac{1}{2\epsilon}\int_0^{t} (L_{(u+\epsilon)\wedge t }^0(X)-L_u^0(X))du=  \frac{1}{2\epsilon} \int_0^{t} \left(\int^{(u+\epsilon)\wedge t }_u dL_s^0(X)\right)du.
$$
Next,  Fubini's theorem gives
\begin{eqnarray*}
\frac{1}{2\epsilon}\int_0^{t} (L_{(u+\epsilon)\wedge t }^0(X)-L_u^0(X)) du   & = & \frac{1}{2}\int_0^{t}  \left(\frac{1}{\epsilon} \int_{(s-\epsilon)^+}^{s }du \right) dL_s^0(X),  \\
 & = &  \frac{1}{2}\int_0^{t}  \frac{s \wedge \epsilon }{\epsilon} dL_s^0(X). 
\end{eqnarray*}
Since $ L^0_t(X)= \int_0^t dL_s^0 (X)$,   we obtain
$$
\left|  \int_0^{t} (L_{(u+\epsilon)\wedge t }^0(X)-L_u^0(X)) du    - L^0_t(X) \right|   \leqslant  \int_0^T \left| \frac{s \wedge \epsilon }{\epsilon}  -1 \right| dL_s^0 (X), \quad 0\leqslant t \leqslant T. 
$$
Observe that $ \left|  \frac{s \wedge \epsilon }{\epsilon} - 1\right| $ is bounded by 2 and tends to 0 for all $s\in ]0,  T]$ as $\epsilon \to 0$. Consequently,   Lebesgue's theorem implies  
 $$ \lim_{\epsilon\to 0 } \frac{1}{2\epsilon}\int_0^{t} (L_{(u+\epsilon)\wedge t }^0(X)-L_u^0(X)) du    -\frac{1}{2}L^0_t(X) = 0, \quad a.s. ,   $$
 uniformly on $t\in[0,  T]$.

   \textbf{1.c.  Study of $\frac{1}{\epsilon}\int_0^t \left( \int_u^{(u+\epsilon)\wedge t }\indi_{\{X_s > 0,  X_u \leqslant 0\} }dX_s \right) du$.}\\
 Obvioulsy, 
\begin{eqnarray*}
 \frac{1}{\epsilon}\int_0^{t} \left( \int_u^{(u+\epsilon)\wedge t }\indi_{\{X_s > 0,   X_u \leqslant 0\} }dX_s\right) du  
 & = & \frac{1}{\epsilon}\int_0^{t} \left( \int_u^{(u+\epsilon)\wedge t }\indi_{\{X_s > 0,   X_u \leqslant 0\} }dV_s\right) du  \\
 &  & +  \frac{1}{\epsilon}\int_0^{t} \left( \int_u^{(u+\epsilon)\wedge t }\indi_{\{X_s > 0,   X_u \leqslant 0\} }dM_s\right) du.
\end{eqnarray*}
Since $V$ has bounded variation, we can proceed as in step \textbf{1.b}. above.
\begin{eqnarray*}
 \left| \frac{1}{\epsilon}\int_0^{t} \left( \int_u^{(u+\epsilon)\wedge t }\indi_{\{X_s > 0,   X_u \leqslant 0\} }dV_s\right) du \right| 
 &=& \left|\int_0^{t} \indi_{\{X_s > 0\} } \left(   \frac{1}{\epsilon} \int^s_{(s-\epsilon)^+}\indi_{\{X_u \leqslant 0\} }du  \right) dV_s \right|, \\
& \leqslant& \int_0^{T} \indi_{\{X_s > 0\} } \left|   \frac{1}{\epsilon} \int^s_{(s-\epsilon)^+}\indi_{\{X_u \leqslant 0\} }du  \right| d |V|_s .
\end{eqnarray*}
Let $s \in [0,T]$ so that $X_s >0$. Since $t\to X_t$ is continuous, $ \frac{1}{\epsilon} \int^s_{(s-\epsilon)^+}\indi_{\{X_u \leqslant 0\} }du=0$ as soon as $\epsilon$ is small enough. Observing that 
$ \left|  \frac{1}{\epsilon} \int^s_{(s-\epsilon)^+}\indi_{\{X_u \leqslant 0\} }du \right|$ is bounded by 1 allows to apply Lebesgue's convergence theorem:
$$\lim_{\epsilon \to 0} \frac{1}{\epsilon}\int_0^{t} \left( \int_u^{(u+\epsilon)\wedge t }\indi_{\{X_s > 0,   X_u \leqslant 0\} }dV_s\right) du =0,$$
a.s, uniformly with respect to $t\in [0,T]$.
 
  Next we claim that 
\begin{equation}
\label{article10eq1}
\lim_{\epsilon \to 0} \textrm{ (ucp) } \frac{1}{\epsilon}\int_0^{t} \left( \int_u^{(u+\epsilon)\wedge t }\indi_{\{X_s > 0,   X_u \leqslant 0\} }dM_s\right) du=0.
\end{equation}

  We may mimic the former analysis which involves $dV_s$, since we have at our disposal a stochastic version of Fubini's theorem (c.f. Lemma \ref{articlefubini}). More precisely,
$$\frac{1}{\epsilon}\int_0^t \left( \int_u^{(u+\epsilon)\wedge t }\indi_{\{X_s > 0,  X_u \leqslant 0\} }dM_s \right) du=
 \int_0^{t} g(\epsilon,   s ) dM_s,  $$
with
$$g(\epsilon,   s )=\indi_{\{ X_s > 0\} }  \frac{1}{\epsilon}\int_{(s-\epsilon)^+}^{s}\indi_{\{  X_u \leqslant 0\} } du . $$
Doob's inequality comes to:
$$
E\left( \sup_{t\in [0,  T]} \left(\int_0^{t} g(\epsilon,   s ) dM_s\right)^2 \right)  \leqslant   4 E\left( \int_0^{T} |g(\epsilon,   s )|^2 d<M>_s \right).  
$$
It is now clear that (\ref{article10eq1}) holds.

 \textbf{1.d. Study of $\frac{1}{\epsilon}\int_0^t \left( \int_u^{(u+\epsilon)\wedge t }\indi_{\{X_s \leqslant 0,  X_u > 0\} }dX_s \right) du$}\\
Although $\frac{1}{\epsilon}\int_0^t \left( \int_u^{(u+\epsilon)\wedge t }\indi_{\{X_s \leqslant 0,  X_u > 0\} }dX_s \right) du$ is quite similar to \\
$\frac{1}{\epsilon}\int_0^{t} \left( \int_u^{(u+\epsilon)\wedge t }\indi_{\{X_s > 0,   X_u \leqslant 0\} }dX_s\right) du$, the difference between large and strict inequality forces us to study it independently. Similary to point \textbf{1.c}, we have:
\begin{eqnarray*}
 \frac{1}{\epsilon}\int_0^{t} \left( \int_u^{(u+\epsilon)\wedge t }\indi_{\{X_s \leqslant 0,  X_u > 0\} }dX_s\right) du  
 & = & \frac{1}{\epsilon}\int_0^{t} \left( \int_u^{(u+\epsilon)\wedge t }\indi_{\{X_s \leqslant 0,  X_u > 0\} }dV_s\right) du  \\
 & + &   \frac{1}{\epsilon}\int_0^{t} \left( \int_u^{(u+\epsilon)\wedge t }\indi_{\{X_s \leqslant 0,  X_u > 0\} }dM_s\right) du.
\end{eqnarray*}
For the first term, we have
$$
 \left| \frac{1}{\epsilon}\int_0^{t} \left( \int_u^{(u+\epsilon)\wedge t }\indi_{\{X_s \leqslant 0,  X_u > 0\} }dV_s\right) du \right| \leqslant
 \int_0^{T} \indi_{\{ X_u > 0\} }   \left| \frac{1}{\epsilon} \int_u^{(u+\epsilon)\wedge t }\indi_{\{X_s \leqslant 0\} }dV_s \right|  du,
 $$
for $0\leqslant t \leqslant T$. Let $u \in [0,T]$ so that $X_u >0$. Using the continuity of $X$, it is clear that $ \frac{1}{\epsilon} \int_u^{(u+\epsilon)\wedge t }\indi_{\{X_s \leqslant 0\} }dV_s=0$ as soon as $\epsilon$ is small enough. Moreover, this term     is bounded by the total variation of $V$.  Thus,  Lebesgue's convergence theorem gives
$$\lim_{\epsilon \to 0}  \frac{1}{\epsilon}\int_0^{T} \left( \int_u^{(u+\epsilon)\wedge t }\indi_{\{X_s \leqslant 0,  X_u > 0\} }dV_s\right) du =0,$$
a.s, uniformly with respect to $t\in [0,T]$.

  For the second term, we proceed as in  point \textbf{1.c}. Fubini's stochastic theorem gives 
$$\frac{1}{\epsilon}\int_0^t \left( \int_u^{(u+\epsilon)\wedge t } \indi_{\{X_s \leqslant 0,  X_u > 0\} }dM_s \right) du= \int_0^{t} \indi_{\{X_s \leqslant 0\} } \tilde{g}(\epsilon,   s ) dM_s,  $$
with
$$\tilde{g}(\epsilon,   s )= \frac{1}{\epsilon}\int_{(s-\epsilon)^+}^{s} \indi_{\{ X_u > 0\} } du . $$
Doob's inequality comes to:
$$
E\left( \sup_{t\in [0,  T]} \left(\int_0^{t} \indi_{\{X_s \leqslant 0\} } \tilde{g}(\epsilon,   s ) dM_s\right)^2 \right)  \leqslant   4 E\left( \int_0^{T}\indi_{\{X_s \leqslant 0\} } | \tilde{g}(\epsilon,   s )|^2 d<M>_s \right).  
$$
We claim that  $E\left( \int_0^{T} \indi_{\{X_s \leqslant 0\} } | \tilde{g}(\epsilon,   s )|^2 d<M>_s \right)$ tends to 0 as $\epsilon \to 0$. Indeed we decompose this term as a sum of two terms:
\begin{eqnarray*}
E\left( \int_0^{T} |\tilde{g}(\epsilon,   s )|^2 d<M>_s \right) &= &E\left( \int_0^{T} \indi_{\{X_s < 0\} } |\tilde{g}(\epsilon, s )|^2 d<M>_s \right)\\
&& + E\left( \int_0^{T} \indi_{\{X_s = 0\} } |\tilde{g}(\epsilon,   s )|^2 d<M>_s \right).
\end{eqnarray*}
Since $X_s < 0$ in the first term, $\tilde{g}(\epsilon, s ) = 0$ as soon as $\epsilon$ is small enough and, by Lebesgue's convergence theorem, this term converges to 0 as $\epsilon \to 0$. Since $\tilde{g}(\epsilon,   s) \leqslant 1$, the second term may be bounded by $E\left( \int_0^{T} \indi_{\{X_s = 0\} } d<M>_s\right)$. By the occupation times formula, 
$$ \int_0^{T} \indi_{\{X_s = 0\} } d<M>_s =  \int_{\R} \indi_{\{x = 0\} } L_T^x dx = 0.$$
 Thus, $ \sup_{t\in [0,  T]} \frac{1}{\epsilon}\int_0^{t} \left( \int_u^{(u+\epsilon)\wedge t } \indi_{\{X_s \leqslant 0,  X_u > 0\} }dM_s\right) du$ tends to 0 in $L^2(\Omega)$ as $\epsilon \to 0$.
 \qed

\begin{rem}
By definition, $I^3_\epsilon(t) $ is a sum of two terms. We have not been able to prove that each term converges. This explains that they have been gathered. Let us explain where comes the difficulty. Let us start with (\ref{article10eq2}) and let use Fubini's theorem in the integral which contains the local time. We get
$$\int_0^{t} (L_{(u+\epsilon)\wedge t }^0(X)-L_u^0(X)) \indi_{\{ X_u \leqslant 0\} } du=\int_0^{t}  \left(\frac{1}{\epsilon} \int_{(s-\epsilon)^+}^{s } \indi_{\{ X_u \leqslant 0\} }du \right) dL_s^0(X).$$
We recall that $ \left(\frac{1}{\epsilon} \int_{(s-\epsilon)^+}^{s } \indi_{\{ X_u \leqslant 0\} }du \right) $ tends to $\indi_{\{ X_s \leqslant 0\} }$ almost surely with respect to  Lebesgue measure. This result cannot guarantee that $ \left(\frac{1}{\epsilon} \int_{(s-\epsilon)^+}^{s } \indi_{\{ X_u \leqslant 0\} }du \right) $ converges $dL_s^0(X)$-everywhere, since the random measure $dL_s^0(X)$ is singular with respect to Lebesgue measure.  
\end{rem}

   \textbf{2. Proof of point \textit{ii)}.}  Let us now suppose that $X$ is a diffusion which satisfies (\ref{article060}) and (\ref{article060b}) .  We will use time reversal.  \\
First, we decompose $I^4_{\epsilon}(t)$ as
\begin{equation}
\label{article10eq3}
I^4_{\epsilon}(t) = \frac{1}{\epsilon}\int_0^{t-\epsilon} X_u^- \indi_{ \{X_{u+\epsilon} > 0\} } du + \frac{1}{\epsilon}\int_0^{t-\epsilon} X_u^+ \indi_{ \{X_{u+\epsilon} \leqslant 0\} } du  +  R^{4,1}_{\epsilon}(t),
\end{equation}
with
$$  R^{4,1}_{\epsilon}(t)  =    \frac{1}{\epsilon}\int_{t-\epsilon}^{t} X_u^- \indi_{ \{X_{ t} > 0\} } du +  \frac{1}{\epsilon}\int_{t-\epsilon}^{t} X_u^+ \indi_{ \{X_{ t} \leqslant 0\} } du.$$
Since $  \frac{1}{\epsilon}\int_{t-\epsilon}^{t} X_u^- \indi_{ \{X_{ t} > 0\} } du =   \frac{1}{\epsilon}\int_{t-\epsilon}^{t}( X_u^- - X_t^-)\indi_{ \{X_{ t} > 0\} } du$, Lemma \ref{articlerestezero} applies and the first term of $R^{4,1}_{\epsilon}(t) $ converges a.s. to 0. The second term in the definition $  R^{4,1}_{\epsilon}(t) $ above converges to 0 by the same way. Thus, $R^{4,1}_{\epsilon}(t) $ converges a.s. to 0 as $\epsilon \to 0$, uniformly on $[0,T]$.

Next, we make the change  of variable $s=T-u-\epsilon$ in the two integrals in (\ref{article10eq3}):
$$
I^4_{\epsilon}(t) =  \frac{1}{\epsilon}\int_{T-t}^{T-\epsilon} X_{T-s-\epsilon}^- \indi_{ \{X_{T-s} > 0\} } ds +  \frac{1}{\epsilon}\int_{T-t}^{T-\epsilon} X_{T-s-\epsilon}^+ \indi_{ \{X_{T-s} \leqslant 0\} } ds  +  R^{4,1}_{\epsilon}(t).
$$
We recall that $\widetilde{X}$ is defined by (\ref{articlereverse}). By Theorem 5.1 of \cite{24}, $(\widetilde{X}_u)_{u\in [0,T]}$ is a semimartingale. We obtain
\begin{eqnarray}
\label{article10eq4}
I^4_{\epsilon}(t) &= & \frac{1}{\epsilon}\int_{T-t}^{T-\epsilon} \widetilde{X}_{s+\epsilon}^- \indi_{ \{\widetilde{X}_{s} > 0\} } ds +   \frac{1}{\epsilon}\int_{T-t}^{T-\epsilon} \widetilde{X}_{s+\epsilon}^+ \indi_{ \{\widetilde{X}_{s} \leqslant 0\} } ds  +  R^{4,1}_{\epsilon}(t), \nonumber \\
& =&
 \left[  \frac{1}{\epsilon}\int_{T-t}^{T} \widetilde{X}_{(s+\epsilon)\wedge T}^- \indi_{ \{\widetilde{X}_{s} > 0\} } ds
  +  \frac{1}{\epsilon}\int_{T-t}^{T} \widetilde{X}_{(s+\epsilon)\wedge T}^+ \indi_{ \{\widetilde{X}_{s} \leqslant 0\} } ds \right]  \label{article10eq3b} \\
  &&
  + R^{4,2}_{\epsilon}(t), \nonumber 
  \end{eqnarray}
with
$$R^{4,2}_{\epsilon}(t) = R^{4,1}_{\epsilon}(t) - \frac{1}{\epsilon}\int_{T-\epsilon}^{T} \widetilde{X}_{T}^- \indi_{ \{\widetilde{X}_{s} > 0\} } ds +  \frac{1}{\epsilon}\int_{T-\epsilon}^{T} \widetilde{X}_{T}^+ \indi_{ \{\widetilde{X}_{s} \leqslant 0\} } ds  .$$
It is clear that $  R^{4,2}_{\epsilon}(t)$ converges a.s to 0 as $\epsilon \to 0$, uniformly on $[0,T]$. 

  Since $\widetilde{X}$ is a semimartingale, point $i)$ of Theorem \ref{article10} may be applied and  the term in bracket in (\ref{article10eq3b}) converges. Moreover, using (\ref{article065b}), we get:
$$ \lim_{\epsilon\to 0}  \textrm{ (ucp) } I^4_{\epsilon}(t) = \frac{1}{2} ( L^0_T(\widetilde{X}) - L^0_{T-t}(\widetilde{X}) )=  \frac{1}{2}  L^0_t(X).
$$
 \qed

\section{Proof of Theorem \ref{article11}  }\label{artsec4}
\setcounter{equation}{0}

   In this section, $X$ will be a standard Brownian motion, and we fix  $\delta\in ]0, \frac{1}{2}[$ (c.f. Lemma \ref{article118}). Since $(-X)$ is a standard Brownian motion, then
$$\lim_{\epsilon \to 0}   I^{3,1}_\epsilon(t) = \lim_{\epsilon \to 0}   I^{3,2}_\epsilon(t),  \qquad
 \lim_{\epsilon \to 0}   I^{4,1}_\epsilon(t) = \lim_{\epsilon \to 0}   I^{4,2}_\epsilon(t) .$$
 Consequently,  we have to prove three distinct properties:
\begin{enumerate}
  \item the (ucp) convergence of $ I^{4,2}_\epsilon(t) $ to $\frac{1}{4}L_t^0(X)$ (see point \textbf{2}.),
  \item the (ucp) convergence of $ I^{3,2}_\epsilon(t) $ to $\frac{1}{4}L_t^0(X)$ (see point \textbf{3}.),
  \item the rate of decay of $ I^{4,i}_\epsilon(t) -\frac{1}{4}L_t^0(X)$ as $\epsilon \to 0$, $i=1,2$ (see point \textbf{2}.).
\end{enumerate}
%
   \textbf{1.}  First, we briefly show that $r^3_\epsilon(t)$ and $r^4_\epsilon(t)$ (defined by (\ref{article09brestec}) and (\ref{article09breste})) tend to 0 as $\epsilon \to 0$, uniformly on $[0,T]$. By the occupation times formula, 
$$\int_0^{t}  \indi_{ \{X_u = 0\} } du = \int_\R \indi_{ \{x = 0\} } L_t^x(X) dx =0.$$
Consequently, $r^3_\epsilon(t)$ vanishes. Next, we decompose $r^4_\epsilon(t)$ and we make the change of variable $v=u+\epsilon$:
$$r^4_\epsilon(t) =\frac{1}{\epsilon}\int_\epsilon^{t} X_{v-\epsilon}^+ \indi_{ \{X_v = 0\} } dv +  \frac{1}{\epsilon}\int_{t-\epsilon}^{t} X_u^+ \indi_{ \{X_t = 0\} } du.$$ 
As previously,  the first integral in $r^4_\epsilon(t) $ vanishes. By Lemma \ref{articlerestezero}, the second term tends a.s. to 0 as $\epsilon \to 0$, uniformly on $[0,T]$, and we have
\begin{equation}
\label{article122}
 | r^4_\epsilon(t) |\leqslant C_\delta \epsilon^\delta.
\end{equation}

   \textbf{2.} Next, we study the convergence of  $ I^{4,2}_\epsilon(t) $ (defined by (\ref{article09b})).\\
Since $ \indi_{\{ X_{(u+\epsilon)\wedge t}<0\}}$ is not $(\shf_u)$-measurable, we "approximate" it by \\
$E(\indi_{\{ X_{(u+\epsilon)\wedge t}<0\}} | \shf_u)$. Let us remark that this term is a function of $X_u$.  Indeed,  
\begin{eqnarray*}
E(\indi_{\{ X_{(u+\epsilon)\wedge t}< 0\}} | \shf_u)  & = & E(\indi_{\{ X_{(u+\epsilon)\wedge t}-X_u< -X_u\}} | X_u),  \\
 & = & \Phi \left(- \frac{X_u}{\sqrt{\epsilon\wedge(t-u)}}\right),
\end{eqnarray*}
where $\Phi$ is the distribution function of the $\shn(0,1)$-law.

  We introduce it in  $ I^{4,2}_\epsilon(t)$. This leads us to consider a new decomposition of $ I^{4,2}_\epsilon(t)$: 
\begin{equation}
\label{article111}
 I^{4,2}_\epsilon(t) =  A^1_{\epsilon}(t)+ A^2_{\epsilon}(t) +    \frac{1}{\epsilon} \int_{0}^t X_{u}^+ \Phi \left( -\frac{X_u}{\sqrt{ \epsilon}}\right)du,  
\end{equation}
where
\begin{eqnarray}
A^1_{\epsilon}(t) & = &\frac{1}{\epsilon} \int_{0}^{t}  X_{u}^+
\left[ \indi_{\{ X_{(u+\epsilon)\wedge t}< 0\}}-   \Phi \left(- \frac{X_u}{\sqrt{\epsilon\wedge(t-u)}}\right) \right] du ,   \label{article111a}\\
A^2_{\epsilon}(t)& = &  \frac{1}{\epsilon} \int_{(t-\epsilon)^+}^{t}  X_{u}^+ \left[  \Phi \left(- \frac{X_u}{\sqrt{t-u}}\right) du -  \Phi \left( -\frac{X_u}{\sqrt{ \epsilon}}\right) \right].  \label{article111b}
\end{eqnarray}
The main term is $ \frac{1}{\epsilon} \int_{0}^t X_{u}^+ \Phi \left( -\frac{X_u}{\sqrt{ \epsilon}}\right)$ since it converges to $\frac{1}{4}L_t^0(X)$ (see point \textbf{2.a}. below).  In point \textbf{2.b},   we will show that $A^1_{\epsilon}(t) + A^2_{\epsilon}(t) $  tends to 0.

   \textbf{2.a.  Study of $\frac{1}{\epsilon} \int_{0}^t X_{u}^+ \Phi \left( -\frac{X_u}{\sqrt{ \epsilon}}\right)du$.} By the occupation times formula,   
$$
\frac{1}{\epsilon} \int_{0}^t X_{u}^+ \Phi \left( -\frac{X_u}{\sqrt{ \epsilon}}\right)du=  \frac{1}{\epsilon} \int_{\R}  x^+ \Phi \left( -\frac{x}{\sqrt{\epsilon}}\right)L^x_t(X)    dx. 
$$
We make the change of variable $y\sqrt{\epsilon}= x$ and we get 
$$\frac{1}{\epsilon} \int_{0}^t X_{u}^+ \Phi \left( -\frac{X_u}{\sqrt{ \epsilon}}\right) du=   \int_0^{\infty}   y \Phi ( -y  ) L^{y\sqrt{\epsilon}}_t(X)    dy.$$
  Since $ \int_0^{\infty}   y \Phi ( -y  ) dy =\frac{1}{4}$,   we have
$$\frac{1}{\epsilon} \int_{0}^t X_{u}^+ \Phi \left( -\frac{X_u}{\sqrt{ \epsilon}}\right)du - \frac{1}{4}L^0_t(X) =  \int_0^{\infty} y \Phi(-y)\left(L_t^{\sqrt{\epsilon} y}(X) - L^0_t (X) \right) dy. $$
By the H\"older property of the Brownian local time  (c.f.  Lemma \ref{article118}),    
\begin{equation}
\label{article121}
\left| \frac{1}{\epsilon} \int_{0}^t X_{u}^+ \Phi \left( -\frac{X_u}{\sqrt{ \epsilon}}\right)- \frac{1}{4}L^0_t(X) \right| \leqslant \left( K_{\delta}\int_0^{\infty} y^{\delta+1}  \Phi(-y)dy\right) \epsilon^{\frac{\delta}{2}} . 
\end{equation}
Therefore $\frac{1}{\epsilon} \int_{0}^t X_{u}^+ \Phi \left( -\frac{X_u}{\sqrt{ \epsilon}}\right)$ converge to $\frac{1}{4}L^0_t(X)$ a.s,   uniformly on $[0,  T]$.

   \textbf{2.b.  Study of $A^1_{\epsilon}(t) + A^2_{\epsilon}(t) $.} First, we modify $A^1_{\epsilon}(t)$ as follows:
\begin{eqnarray*}
A^1_{\epsilon}(t) & = &\frac{1}{\epsilon} \int_{0}^{(t-\epsilon)^+}  X_{u}^+
\left[ \indi_{\{ X_{u+\epsilon}< 0\}} -\Phi \left( -\frac{X_u}{\sqrt{\epsilon}}\right)\right] du\\
&& + \frac{1}{\epsilon} \int_{(t-\epsilon)^+}^{t}  X_{u}^+
\left[ \indi_{\{ X_{t}< 0\}} -\Phi \left( -\frac{X_u}{\sqrt{t-u}}\right)\right] du.
\end{eqnarray*}
Secondly, we decompose $ \indi_{\{ X_{t}< 0\}} -\Phi \left( -\frac{X_u}{\sqrt{t-u}}\right) $ as:
\begin{eqnarray*} 
&&\indi_{\{ X_{t}< 0\}} - \Phi\left( \frac{-X_t}{\sqrt{u+\epsilon -t}}\right) 
+  \left[  \Phi\left( \frac{-X_u}{\sqrt{\epsilon}}\right)  -\Phi \left(  \frac{-X_u}{\sqrt{t-u}}\right)\right] \\
&& + \left[  \Phi\left( \frac{-X_t}{\sqrt{u+\epsilon -t}}\right)  - \Phi\left(  \frac{-X_u}{\sqrt{\epsilon}}\right) \right]. 
\end{eqnarray*}
Multiplying by $X_{u}^+$ and integrating over $[(t-\epsilon)^+,t]$ give rise to four integrals, but we remark that the third one is equal to $-A^2_{\epsilon}(t)$. Then,
\begin{equation}
\label{article112}
A^1_{\epsilon}(t) + A^2_{\epsilon}(t) = R^1_{\epsilon}(t) + D^2_{\epsilon}(t)+ D^3_{\epsilon}(t) + R^4_{\epsilon}(t) ,  
\end{equation}
where
\begin{eqnarray*}
R^1_{\epsilon}(t) & = &\frac{1}{\epsilon} \int_{0}^{(t-\epsilon)^+}  X_{u}^+
\left[ \indi_{\{ X_{u+\epsilon}< 0\}} -\Phi \left( -\frac{X_u}{\sqrt{\epsilon}}\right)\right] du,  
\\
D^2_{\epsilon}(t) & = &   \indi_{\{ X_{t}< 0\}}  \frac{1}{\epsilon} \int_{(t-\epsilon)^+}^{t}  X_{u}^+du,   \\
D^3_{\epsilon}(t) & = & - \frac{1}{\epsilon} \int_{(t-\epsilon)^+}^{t}  X_{u}^+  \Phi\left( \frac{-X_t}{\sqrt{u+\epsilon -t}}\right)  du,\\
R^4_{\epsilon}(t) &= & \frac{1}{\epsilon} \int_{(t-\epsilon)^+}^{t}  X_{u}^+
\left[  \Phi\left( \frac{-X_t}{\sqrt{u+\epsilon -t}}\right)  - \Phi\left(  \frac{-X_u}{\sqrt{\epsilon}}\right) \right] du.
\end{eqnarray*}
For $u$ fixed,   It\^o's formula applied to $\phi (x,  s) =\Phi\left( -\frac{x}{\sqrt{u+\epsilon -s}}\right)$,\\
  $(\phi (x,  s)  \in \shc^2( ]0,+\infty[, [u, u+\epsilon[ )$,  comes to:
$$
\phi (X_v,  v) -\phi (X_u,  u)=
-\int_{u}^{v}  \frac{\exp \left( -\frac{X_s^2}{2(u+\epsilon -s)}\right)}{\sqrt{2\pi (u+\epsilon -s) }}dX_s, \quad \forall v\in [u,   u+\epsilon[. 
$$
Since the Lebesgue measure of $\{ u; X_{u+\epsilon} =0\}$ vanishes, this formula, used for $v= u+\epsilon$ if $u\in [0,   (t-\epsilon)^+]$, and for $v=t$ if $u\in ](t-\epsilon)^+,   t]$, leads to
\begin{eqnarray}
 \indi_{\{ X_{u+\epsilon}< 0\}} -\Phi \left( -\frac{X_u}{\sqrt{\epsilon}}\right) &=&
-\int_{u}^{u+\epsilon}  \frac{\exp \left( -\frac{X_s^2}{2(u+\epsilon -s)}\right)}{\sqrt{2\pi (u+\epsilon -s) }}dX_s,  \label{article112a} \\
\Phi\left(- \frac{X_t}{\sqrt{u+\epsilon -t}}\right)  - \Phi\left(- \frac{X_u}{\sqrt{\epsilon}}\right) & = &
-\int_{u}^{t}  \frac{\exp \left( -\frac{X_s^2}{2(u+\epsilon -s)}\right)}{\sqrt{2\pi (u+\epsilon -s) }}dX_s.  \label{article112b}
\end{eqnarray}
Using (\ref{article112a}) and   (\ref{article112b}) permits to obtain a new expression of $A^1_{\epsilon}(t) + A^2_{\epsilon}(t) $:
 \begin{equation}
\label{article113}
A^1_{\epsilon}(t) + A^2_{\epsilon}(t) = D^1_{\epsilon}(t) + D^2_{\epsilon}(t)+ D^3_{\epsilon}(t),  
\end{equation}
where
$$ D^1_{\epsilon}(t) =  - \frac{1}{\epsilon} \int_0^{t}  X_{u}^+\left( \int_{u}^{(u+\epsilon)\wedge t}  \frac{\exp \left( -\frac{X_s^2}{2(u+\epsilon -s)}\right)}{\sqrt{2\pi (u+\epsilon -s) }}dX_s\right) du.
$$
In the rest of the paragraph, we will prove that each term $ D^i_{\epsilon}(t) $ tends to 0,   uniformly on the  compact set.  Let us observe that, in $ D^2_{\epsilon}(t)$ and $  D^3_{\epsilon}(t)$,  the lenghts of the intervals of integration are  smaller than $\epsilon$. It will be shown that   the convergence of these terms holds almost surely and is a consequence of the continuity of $X$.  As for the first term $D^1_{\epsilon}(t)$,  we observe that it is a martingale (via Fubini's theorem) and we  use Doob's inequality to get the convergence in the $L^2$ sense.  

   \textbf{Convergence of $D^2_{\epsilon}(t)$ to 0.} Since $\indi_{\{ X_{t}\leqslant 0\}}  X^+_t = 0$,   we have
$$D^2_{\epsilon}(t)  =   \left(\frac{1}{\epsilon} \int_{(t-\epsilon)^+}^{t} (X_{u}^+ -X_{t}^+) \indi_{\{ X_{t}< 0\}} du\right). $$
By Lemma \ref{articlerestezero},  $ D^2_{\epsilon}(t)$ tends a.s.  to $0$,   uniformly on $[0,  T]$ and we have
\begin{equation}
\sup_{t\in[0,T]}\left| D^2_{\epsilon}(t) \right|    \leqslant  | C_{\delta}| \epsilon^{\delta}.  \label{article115}
\end{equation}

   \textbf{Convergence of $D^3_{\epsilon}(t)$ to 0.} Let us introduce $X^+_t $  in $D^3_{\epsilon}(t)$, we get
   $$\begin{array}{rl}
  D^3_{\epsilon}(t) = &   - \frac{1}{\epsilon} \int_{(t-\epsilon)^+}^{t} (X_t^+ - X_{u}^+)  \Phi\left( \frac{-X_t}{\sqrt{u+\epsilon -t}}\right)  du  \\
  &  +  \frac{1}{\epsilon} \int_{(t-\epsilon)^+}^{t}  X_{t}^+  \Phi\left( \frac{-X_t}{\sqrt{u+\epsilon -t}}\right)  du.  
\end{array}
 $$
%
In the first term,  by using the  fact that  $\Phi$ is bounded by 1 and Lemma \ref{article118} (H\"older property of $t\to X_t$), we get
$$
\left| \frac{1}{\epsilon} \int_{(t-\epsilon)^+}^{t} (X_t^+ - X_{u}^+)  \Phi\left( \frac{-X_t}{\sqrt{u+\epsilon -t}}\right)  du \right| \leqslant   \frac{|C_\delta |}{\epsilon} \int_{(t-\epsilon)^+}^{t}  |t-u|^\delta du
 \leqslant  |C_\delta | \epsilon^\delta. 
$$
Let us consider the second term in $D^3_{\epsilon}(t)$. Since   $ \Phi(x) \leqslant C e^{\frac{-x^2}{4}}$ for any $x \leqslant 0$, we obtain:
$$\left|  \frac{1}{\epsilon} \int_{(t-\epsilon)^+}^{t}  X_t^+\Phi \left( \frac{-X_t}{ \sqrt{u+\epsilon -t } } \right) du \right|\leqslant  \frac{C}{\epsilon} \int_{(t-\epsilon)^+}^{t}  X_{t}^+ e^{-\frac{X_t^2}{4(u+\epsilon -t)} } du. $$
The function defined on $\R^+$ by  $x\to x  e^{-\frac{x^2}{4(u+\epsilon -t)} } $ is non-negative and bounded by  $\sqrt{2 (u+\epsilon -t )}e^{-\frac{1}{2}}$.  Consequently,   
$$
\left|  \frac{1}{\epsilon} \int_{(t-\epsilon)^+}^{t}  X_t^+\Phi \left( \frac{-X_t}{ \sqrt{u+\epsilon -t } } \right) du \right|\leqslant C\frac{1}{\epsilon} \int_{(t-\epsilon)^+}^{t} \sqrt{(u+\epsilon -t )} du \leqslant C\sqrt{\epsilon}. 
$$
Thus, $D^3_{\epsilon}(t)$ converge a.s.  to 0 uniformly on $ [0,  T]$ and
\begin{equation}
\label{article116}
\sup_{t\in [0,  T]} |D^3_{\epsilon}(t)| \leqslant  |C_\delta | \epsilon^\delta + C\sqrt{\epsilon}. 
\end{equation}

   \textbf{Convergence of $D^1_{\epsilon}(t)$ to 0.}  In order to use Fubini's  theorem (i.e. Theorem \ref{articlefubini}),  first we prove that
$$ \int_0^{t}  \left( \int _u^{(u+\epsilon)\wedge t}  E \left[ H^2_{\epsilon}(u,  s) \right] ds\right) du <\infty,  $$
where
$$H_{\epsilon}(u,  s)= \frac{X_u^+}{\epsilon \sqrt{2\pi (u+\epsilon-s) }}\exp \left( -\frac{X_s^2}{2(u+\epsilon-s)}\right)   . $$
By Cauchy-Schwarz inequality,   we have
$$
 E \left[  H^2_{\epsilon}(u,  s) \right]   \leqslant  \frac{1}{2\pi (u+\epsilon-s) \epsilon^2} \sqrt{ E[(X_u^+ )^4]  E \left[ \exp \left( -\frac{2 X_s^2}{(u+\epsilon-s)}\right) \right] }. $$
Since $X_r$ is a centered Gaussian random variable with variance $r$,   we can calculate explicitly  the two expectations which appear in the right-hand side of the previous inequality:
$$E[(X_u^+ )^4]  = \frac{3}{2}u^2, \quad E \left[ \exp \left( -\frac{2 X_s^2}{(u+\epsilon-s)}\right) \right] = \sqrt{\frac{u+\epsilon-s}{3s + u + \epsilon}}.$$
Since $s \geqslant u \Rightarrow 3s + u + \epsilon > 4 u$, reporting in the inequality gives
$$
 E \left[  H^2_{\epsilon}(u,  s) \right]   \leqslant  \frac{C}{\epsilon^2} \left( \frac{u}{u+\epsilon -s}\right)^{\frac{3}{4}},   \quad \forall u \in [0,  t],   s \in]u,   u+\epsilon[. $$
By an easy calculation, it may be proved that $ \int_0^{t} \left( \int _u^{(u+\epsilon)\wedge t}  \left( \frac{u}{u+\epsilon -s}\right)^{\frac{3}{4}}ds \right )du < \infty$. Consequently,    Fubini's stochastic theorem may be applied, we have:
$$D^1_{\epsilon}(t)  = - \int _0^{t} \left( \int_{(s-\epsilon)^+}^{s}   \frac{ 1 }{\epsilon}\frac{X_u^+}{\sqrt{2\pi (u+\epsilon -s)}}\exp \left( -\frac{X_s^2}{2(u+\epsilon-s)}\right)  du\right) dX_s. $$
Thus,   $(D^1_{\epsilon}(t),   0\leqslant t \leqslant T) $ is a square integrable martingale,   and by Doob's inequality, we get:
$$E\left(\sup_{t\in [0,  T]} ( D^1_{\epsilon}(t)  )^2 \right) \leqslant 4 E\left[( D^1_{\epsilon} (T))^2 \right], $$
and
$$
E\left[( D^1_{\epsilon} (T))^2 \right] =
 E\left( \int _0^{T} \left[ \int_{(s-\epsilon)^+}^{s}  \frac{X_u^+ \exp \left( -\frac{X_s^2}{2(u+\epsilon -s)}\right)}{\epsilon\sqrt{2\pi (u+\epsilon -s)}} du \right]^2 ds  \right) . $$
 %
The Cauchy-Schwarz inequality applied to the term in brackets gives:
$$
E\left[( D^1_{\epsilon} (T))^2 \right]  \leqslant 
  \int _0^{T}   \left[   \int_{(s-\epsilon)^+}^{s} \frac{1}{2\pi \epsilon (u+\epsilon -s)}
   E\left((X_u^+)^2  e^{-\frac{X_s^2}{u+\epsilon -s}}  \right)du \right] ds . 
$$
Let $s > u > 0$, we introduce $A:=  E\left((X_u^+)^2 e^{-\frac{X_s^2}{u+\epsilon -s}}   \right)$. Since $X_s -X_u$ is independent from $X_u$, we have
$$ A=  \frac{1}{2\pi \sqrt{(s-u)u}} \int_0^{\infty} x^2 e^{ -\frac{x^2}{2u}} \int_{\R} e^{ -\frac{(y+ x)^2}{u+\epsilon -s} -\frac{y^2}{2(s-u)} }dy dx.
$$
Since
$$ \frac{(y+x)^2}{u+\epsilon -s} + \frac{y^2}{2(s-u)} = \frac{s-u+\epsilon}{2(s-u)(u+\epsilon -s)}\left( y + \frac{2(s-u) x}{s-u+\epsilon} \right) +  \frac{x^2}{s-u+\epsilon},$$
we obtain
$$
A = \frac{1}{ \sqrt{2\pi}} \sqrt{\frac{u+\epsilon -s}{u(s-u+\epsilon)}}\int_0^{\infty} x^2 e^{ -\frac{(s+u+\epsilon)x^2}{2u(s-u+\epsilon)}} dx=\frac{u(s-u+\epsilon)\sqrt{u+\epsilon -s}}{2 (s+u+\epsilon)^{\frac{3}{2}}} .$$
Since $u \geqslant s-\epsilon \Rightarrow s-u+\epsilon \leqslant 2\epsilon$ and $\frac{1}{s+u+\epsilon } \leqslant\frac{1}{u}$, reporting in the integral gives:
$$
E\left[( D^1_{\epsilon} (T))^2 \right] \leqslant  C
  \int _0^{T}  \left( \int_{u}^{u+\epsilon} \frac{1}{\sqrt{u+\epsilon -s}} ds\right) \frac{ du }{\sqrt{u}} \leqslant C \sqrt{T \epsilon} .   $$
Hence,  
\begin{equation}
\label{article120}
E\left[ \sup_{t\in [0,  T]} \left( D^1_{\epsilon}(t)  \right)^2\right] \leqslant C \sqrt{T\epsilon},   
\end{equation}
and $D^1_{\epsilon}(t)$ converge to 0 in $L^2(\Omega)$,   uniformly on $[0,  T]$.

   \textbf{2.c.} We are now able to prove point \textbf{2}. of Theorem \ref{article11}. It is clear that (\ref{article111}) and (\ref{article113})  imply:
\begin{equation}
\label{article123}
I^{4,2}_\epsilon(t)  - \frac{1}{4}L^0_t(X) = D^1_{\epsilon}(t) + D^2_{\epsilon}(t)+ D^3_{\epsilon}(t) +  \left(  \frac{1}{\epsilon} \int_{0}^t X_{u}^+ \Phi \left( -\frac{X_u}{\sqrt{ \epsilon}}\right) - \frac{1}{4}L^0_t(X) \right).
\end{equation} 
The (ucp) convergence of $I^{4,2}_\epsilon(t)$ to  $ \frac{1}{4}L^0_t(X)$ and
 $$\left\| \sup_{t\in [0,  T]} \left| I^{4,2}_\epsilon(t) -   \frac{1}{4}L^0_t(X) \right| \right\|_{L^2(\Omega)}
  \leqslant  C \epsilon^{\frac{\delta}{2}}, $$
  is a direct consequence of inequalities (\ref{article115}),  (\ref{article116}),   (\ref{article120}) and (\ref{article121}).

   \textbf{3.} Finally, we study the convergence of  $ I^{3,2}_\epsilon(t) $ (defined by (\ref{article09c})).

  We follow the same method as in point \textbf{2}. Since $X_{u+\epsilon}^+$ is not $(\shf_u)$-measurable, we "replace" it by $E( X_{u+\epsilon}^+  | \shf_u)$. This leads us to consider the following  decomposition of $ I^{3,2}_\epsilon(t)$: 
$$ I^{3,2}_\epsilon(t) =   \widetilde{I}^{3,2}_\epsilon(t) +  \overline{I}^{3,2}_\epsilon(t)  +  \Delta^{3,2}_\epsilon(t),$$
where
\begin{eqnarray}
  \widetilde{I}^{3,2}_\epsilon(t) & = &  \frac{1}{\epsilon}\int_0^{t} E( X_{u+\epsilon}^+  | \shf_u) \indi_{ \{X_u < 0\} }  du    ,   \label{article111ac}\\
 \overline{I}^{3,2}_\epsilon(t) & = &  \frac{1}{\epsilon}\int_0^{(t-\epsilon)^+} \left( X_{u+\epsilon}^+- E( X_{u+\epsilon}^+  | \shf_u)\right) \indi_{ \{X_u < 0\} }  du \\
 &&
 +  \frac{1}{\epsilon}\int_{(t-\epsilon)^+}^{t} \left( E( X_{u+\epsilon}^+  | \shf_t) - E( X_{u+\epsilon}^+  | \shf_u)\right) \indi_{ \{X_u < 0\} }  du
  ,  \label{article111bc}\\
  \Delta^{3,2}_\epsilon(t) &=&  \frac{1}{\epsilon}\int_{(t-\epsilon)^+}^{t} \left( X_t^+ -E( X_{u+\epsilon}^+  | \shf_t) \right) \indi_{ \{X_u < 0\} }  du .   \label{article111cc}
\end{eqnarray}
The main term is  $  \widetilde{I}^{3,2}_\epsilon(t) $ since it converges to $\frac{1}{4}L_t^0(X)$ (see point \textbf{3.a}. below).  In point \textbf{3.b} ,   we will show that $ \overline{I}^{3,2}_\epsilon(t) $   tends to 0.  As for $\Delta^{3,2}_\epsilon(t)$, since $X_t^+ -E( X_{u+\epsilon}^+  | \shf_t) = Y_t - Y_{u+\epsilon}^+ $, with $Y_s= E( X_s^+  | \shf_t), s\in [(t-\epsilon)^+,t]$, Lemma \ref{articlerestezero} applies and $\Delta^{3,2}_\epsilon(t)$ converges a.s. to 0 as $\epsilon \to 0$, uniformly on $[0,T]$.

   \textbf{3.a.  Study of $  \widetilde{I}^{3,2}_\epsilon(t) $.} First, let us remark that $E( X_{u+\epsilon}^+  | \shf_u)$ is a function of $X_u$.  Indeed,  
$$
E( X_{u+\epsilon}^+  | \shf_u)  =  \sqrt{\epsilon} g\left( \frac{X_u}{ \sqrt{\epsilon}} \right),
$$
where $g(x)= E( (G + x)^+ )$ and $G$ is a Gaussian random variable with $\shn(0,1)$-law.

By the occupation times formula, we have:  
$$
  \widetilde{I}^{3,2}_\epsilon(t) =  \frac{1}{\sqrt{\epsilon}}\int_{-\infty}^0 g\left( \frac{x}{ \sqrt{\epsilon}} \right) L_t^x(X) dx. 
$$
The change of variable $y\sqrt{\epsilon}= x$ gives
$$  \widetilde{I}^{3,2}_\epsilon(t) =   \int_{-\infty}^0 g(y) L_t^{ \sqrt{\epsilon}y}(X) dy.$$
Since 
$$
\begin{array}{rl}
\int_{-\infty}^0 g(y) dy  =& E\left( \int_{-\infty}^0 ( G + y)^+ dy \right)  = E\left(\indi_{\{ G >0 \}  }\int_{-G}^0 (G + y) dy  \right),\\
= &\frac{1}{2} E( (G^+)^2) = \frac{1}{4},
\end{array}$$
 %
  we have
$$  \widetilde{I}^{3,2}_\epsilon(t)  - \frac{1}{4}L_t^0(X) =   \int_{-\infty}^0 g(y) ( L_t^{ \sqrt{\epsilon}y}(X) - L_t^0(X)) dy.$$
Consequently,  the H\"older property of the Brownian local time  (c.f.  Lemma \ref{article118}), implies that:   
$$ \left| \widetilde{I}^{3,2}_\epsilon(t)  - \frac{1}{4}L_t^0(X) \right| \leqslant \epsilon ^{\frac{\delta}{2}} K_\delta E\left( \int_{-\infty}^0    ( G + y)^+   y^\delta dy \right) \leqslant C\epsilon ^{\frac{\delta}{2}}  .$$
Therefore $\widetilde{I}^{3,2}_\epsilon(t) $ converge to $\frac{1}{4}L^0_t(X)$ a.s,   uniformly on $[0,  T]$.

   \textbf{3.b  Study of $  \overline{I}^{3,2}_\epsilon(t) $.} Our goal is to show that $  \overline{I}^{3,2}_\epsilon(t) $ is a martingale, in order to apply Doob's inequality.  To begin with, we write $ \left( X_{u+\epsilon}^+- E( X_{u+\epsilon}^+  | \shf_u)\right)$ and   $ \left( E( X_{u+\epsilon}^+  | \shf_t) - E( X_{u+\epsilon}^+  | \shf_u)\right) $ as stochastic integrals using Ito's formula.

  Let $u, \epsilon$ be fixed. We define $M_s =E( X_{u+\epsilon}^+  | \shf_s), s\in [u,u+\epsilon[$.  $(M_s)_{s\in [u,u+\epsilon[}$ is a martingale and
$$M_s =E( (X_{u+\epsilon}-X_s + X_s)^+  | \shf_s) =  f( s, X_s),$$ 
where 
$$
f(s, y) = \int_{\R} (x+y)^+ p_{u+\epsilon -s} (x) dx = \int_0^{+\infty} z p_{u+\epsilon -s} (z-y) dz,
$$
and $p_{u+\epsilon -s}(x)= \frac{e^{-\frac{x^2}{2(u+\epsilon -s)}}}{\sqrt{2\pi (u+\epsilon -s) }}$.

Since $\frac{1}{2}  \frac{\partial^2 f }{\partial y^2}(s, y)  = -\frac{\partial f }{\partial s}(s, y)$, applying  It\^o's formula to $f\in \shc^{1,2} ( [u,u+\epsilon[, \R)$ comes to:
\begin{eqnarray}
 f( u+ \epsilon, X_{u+ \epsilon}) - f ( u, X_u) &=&\int_u^{u+ \epsilon} \frac{\partial f }{\partial y} (s, X_s) dX_s, \quad u\in [0, (t-\epsilon)^+],  \label{articlerajout1}\\
  f( t, X_t) - f ( u, X_u) &=&\int_u^t \frac{\partial f }{\partial y} (s, X_s) dX_s, \quad u \in [ (t-\epsilon)^+ , t] . \label{articlerajout2}
\end{eqnarray}
Let us evaluate $\frac{\partial f }{\partial y}(s,y)$. Writing $ f(s,y) =  \int_{-y}^{+\infty} (x+y) p_{u+\epsilon -s} (x) dx$ allows to calculate the $y$-derivative of $f$:
$$
\frac{\partial f }{\partial y}(s,y) =  -  \int_{-y}^{+\infty} p_{u+\epsilon -s}(x)dx =  1- \Phi \left( \frac{-y}{\sqrt{u+\epsilon -s}} \right) ,
$$
where $\Phi$ is the distribution function of the standard Gaussian distribution.

Let us observe that
$$
 f( u+ \epsilon, X_{u+ \epsilon}) - f ( u, X_u) = M_{u+ \epsilon} - M_u = X_{u+\epsilon}^+ - E( X_{u+\epsilon}^+  | \shf_u), $$
 and
 $$
  f( t, X_t) - f ( u, X_u) = M_t - M_u =E( X_{u+\epsilon}^+  | \shf_t) - E( X_{u+\epsilon}^+  | \shf_u).
$$
As a result, reporting   (\ref{articlerajout1}) and  (\ref{articlerajout2}) in  $\overline{I}^{3,2}_\epsilon(t)$ gives:
$$
 \overline{I}^{3,2}_\epsilon(t) =   \frac{1}{\epsilon}\int_0^{t} 
\left[ \int_u^{(u+ \epsilon)\wedge t} ( 1- \Phi \left( \frac{-X_s}{\sqrt{u+\epsilon -s}} \right) )dX_s\right] 
\indi_{ \{X_u < 0\} }   du.
$$
Since $x\to 1- \Phi(x) $ is uniformly bounded by $1$,  Fubini's  theorem (i.e. Theorem \ref{articlefubini}) may be applied, we have
 $$ \overline{I}^{3,2}_\epsilon(t) = \int_0^{t} 
\left[  \frac{1}{\epsilon} \int_{(s-\epsilon)^+}^s (1- \Phi \left( \frac{-X_s}{\sqrt{u+\epsilon -s}} \right)  )\indi_{ \{X_u < 0\} }   du\right] dX_s.$$
Thus,   $( \overline{I}^{3,2}_\epsilon(t),   0\leqslant t \leqslant T) $ is a square integrable martingale.  By Doob's inequality, we get:
$$ E\left( \sup_{t\in[0,T]}  |\overline{I}^{3,2}_\epsilon(t)|^2 \right) \leqslant 4 E\left(  \int_0^{T} 
\left[  \frac{1}{\epsilon} \int_{(s-\epsilon)^+}^s (1- \Phi \left( \frac{-X_s}{\sqrt{u+\epsilon -s}} \right)  )\indi_{ \{X_u < 0\} }   du\right]^2 ds \right).$$
Finally, we will show that the right-hand side of the inequality above converges to 0. Recall that $\int_0^T \indi_{ \{X_s = 0\} }  ds =0$, let us introduce two cases: $X_s <0$ and $X_s > 0$. 
\begin{itemize}
  \item We have:
 $$ \indi_{ \{X_s > 0\} }  \left| \frac{1}{\epsilon} \int_{(s-\epsilon)^+}^s ( 1- \Phi \left( \frac{-X_s}{\sqrt{u+\epsilon -s}} \right)  )\indi_{ \{X_u < 0\} }   du \right| $$
 $$ \leqslant  \indi_{ \{X_s > 0\} } \left(  \frac{1}{\epsilon} \int_{(s-\epsilon)^+}^s \indi_{ \{X_u < 0\} } du\right).$$
We have already proven that the right-hand side of the prior inequality goes to 0, as $\epsilon \to 0$, a.s. and in $L^1(\Omega)$.
  \item if $X_s < 0$, by using $|1-\Phi (\alpha)| \leqslant  C e^{-\frac{\alpha^2}{4}}$ for any $\alpha > 0$, we get
   $$ \left| \frac{1}{\epsilon} \int_{(s-\epsilon)^+}^s ( 1- \Phi \left( \frac{-X_s}{\sqrt{u+\epsilon -s}} \right)  )\indi_{ \{X_u < 0\} }   du \right| \leqslant  \frac{C}{\epsilon} \int_{(s-\epsilon)^+}^s e^{-\frac{X_s^2}{4 (u+\epsilon -s)} } du.$$
The change of variable $v = u+\epsilon -s$ gives:
   $$ \left| \frac{1}{\epsilon} \int_{(s-\epsilon)^+}^s ( 1- \Phi \left( \frac{-X_s}{\sqrt{u+\epsilon -s}} \right)  )\indi_{ \{X_u < 0\} }   du \right| \leqslant  \frac{C}{\epsilon} \int_{(\epsilon-s)^+}^{\epsilon} e^{-\frac{X_s^2}{4v} } dv.$$
Since   $X_s \neq 0$, 
$\frac{1}{\epsilon} \int_{(\epsilon-s)^+}^{\epsilon} e^{-\frac{X_s^2}{4v} } dv$ converges to 0 a.s and in $L^1(\Omega)$, as $\epsilon \to 0$.
\end{itemize}
Then, by Lebesgue Theorem, $ E\left(  \int_0^{T} \left[  \frac{1}{\epsilon} \int_{(s-\epsilon)^+}^s ( 1- \Phi \left( \frac{-X_s}{\sqrt{u+\epsilon -s}} \right)  )\indi_{ \{X_u < 0\} }   du\right]^2 ds \right)$ converges to 0.  \qed

\section{Proofs of  Theorem \ref{article06b}  and Proposition \ref{article06bb} }\label{artsec5}
\setcounter{equation}{0}

   \textbf{1.} In this Section,   $X$ is supposed to be the standard Brownian motion. It is convenient to adopt the convention that $X_s = 0$ if $s< 0$. The proof of  Theorem \ref{article06b} is based on the identity
\begin{equation}
\label{article0613}
J_{\epsilon}(t)-L_t^0(X) = -\left(I^1_{\epsilon}(t)-  \int_0^t \indi_{\{ 0<X_{s}\}} dX_s \right) + \left(   I^2_{\epsilon}(t) -  X_t^+ -\frac{1}{2}L_t^0(X)\right). 
\end{equation}
In step  \textbf{2.} (resp. \textbf{3.}) below, we study the convergence of $I^1_{\epsilon}(t)$ (resp. $I^2_{\epsilon}(t)$). We will use Theorem  \ref{article11} to obtain the convergence of $I^2_{\epsilon}(t)$ to its limit. In step \textbf{4.} we will show Proposition \ref{article06bb}.

   \textbf{2. Study of the convergence of $I^1_{\epsilon}(t)$.} We will use  (\ref{article068}).  Let us remark that $ \widehat{I}^1_{\epsilon}(t)=0$  since $X$ is a martingale. Therefore  (\ref{article068}) reduces to:
\begin{equation}
\label{article068brownien}
 I^1_{\epsilon}(t)-\int_0^t \indi_{\{ 0<X_{s}\}} dX_s = \widetilde{I}^1_{\epsilon}(t) +\Delta_1(t,  \epsilon),  
\end{equation}
where  $\widetilde{I}^1_{\epsilon}(t),   \Delta_1(t,  \epsilon)$ are defined by (\ref{article069}), respectively (\ref{article0610}). As for $ \Delta_1(t,  \epsilon) $,   Lemma \ref{articlerestezero} gives:
\begin{equation}
\label{article0619}
 \sup_{t\in [0,  T]} | \Delta_1(t,  \epsilon) | \leqslant  C_\delta \epsilon^\delta. 
\end{equation}
Let us now deal with $\widetilde{I}^1_{\epsilon}(t)$. Recall that $<X>_u = u$, the inequality  (\ref{article0611}) becomes: 
$$
E\left(  \sup_{0\leqslant t \leqslant T}( \widetilde{I}^1_{\epsilon}(t) )^2 \right)  \leqslant  4 E \left[\int_0^T  \left(\frac{1}{\epsilon} \int_{(u-\epsilon)^+}^{u}  \indi_{\{ 0<X_{s}\}}- \indi_{\{ 0<X_{u}\}} ds \right)^2 du\right].  
$$
By writing $ \left(\int_{(u-\epsilon)^+}^{u}  \indi_{\{ 0<X_{s}\}}- \indi_{\{ 0<X_{u}\}} ds \right)^2$ as 
$$2 \iint_{[(u-\epsilon)^+,u]^2}\indi_{\{s<s' \}}
(\indi_{\{ 0<X_{s}\}}  - \indi_{\{ 0<X_{u}\}}) ( \indi_{\{ 0<X_{s'}\}}  - \indi_{\{ 0<X_{u}\}}) dsds',$$
 we obtain:
$$
\begin{array}{l}
E\left(  \sup_{0\leqslant t \leqslant T}( \widetilde{I}^1_{\epsilon}(t) )^2 \right) \leqslant   \int_0^T \left\{  \iint_{[(u-\epsilon)^+,u]^2} \frac{8\indi_{\{s<s' \}}}{\epsilon^2} \right.    \\
  \left. E\left ( (\indi_{\{ 0<X_{s}\}}  - \indi_{\{ 0<X_{u}\}}) ( \indi_{\{ 0<X_{s'}\}}  - \indi_{\{ 0<X_{u}\}})\right) dsds' \right\} du. 
\end{array}
$$
%
%
The expectation in the previous integral may be computed explicitly:
$$\begin{array}{l}
E\left ( 
(\indi_{\{ 0<X_{s}\}}  - \indi_{\{ 0<X_{u}\}}) ( \indi_{\{ 0<X_{s'}\}}  - \indi_{\{ 0<X_{u}\}})\right)   \\
= P( 0< X_s, 0< X_{s'}) -P( 0< X_s, 0< X_{u}) -P( 0< X_u, 0< X_{s'})  \\
\quad + P( 0< X_u), \\
= \frac{1}{2\pi} \left(f(\frac{u-s'}{s'})+ f(\frac{u-s}{s}) -f(\frac{s'-s}{s}) \right),   
\end{array}$$
%
%
with $f(x)=\textrm{Arctan}(\sqrt{x})$. Consequently, we have to determine the upper bound of
$$
\int_0^T \left\{  \iint_{[(u-\epsilon)^+,u]} \frac{\indi_{\{s<s' \}}}{\epsilon^2}\left(f(\frac{u-s'}{s'})+ f(\frac{u-s}{s}) -f(\frac{s'-s}{s}) \right)dsds' \right\} du.
$$
The integral may be calculated, then bounded via the inequality $f(x)\leqslant \sqrt{x}$ (see \cite{0} for the details of this fastidious calculus). We get 
\begin{equation}
\label{article0615}
E\left(  \sup_{0\leqslant t \leqslant T}( \widetilde{I}^1_{\epsilon}(t) )^2 \right) \leqslant  C \sqrt{T} \sqrt{\epsilon}. 
\end{equation}

   \textbf{3. Study of the convergence of $I^2_{\epsilon}(t)$.} First, we decompose $I^2_{\epsilon}(t)- X_t^+ -\frac{1}{2}L_t^0(X)$ as:
\begin{equation}
\label{article0616}
I^2_{\epsilon}(t)-  X_t^+ -\frac{1}{2}L_t^0(X)= \left( \widetilde{I}^2_{\epsilon}(t) -  X_t^+ -\frac{1}{2}L_t^0(X) \right)+\Delta_2(t,  \epsilon),  
\end{equation}
where
\begin{eqnarray*}
 \widetilde{I}^2_{\epsilon}(t) &=&  \int_0^{t} \frac{1}{\epsilon}(X_{(s+\epsilon)\wedge t}-X_{s}) \indi_{\{ 0<X_{(s+\epsilon)\wedge t} \}} ds ,  \\
 \Delta_2(t,  \epsilon)&=&  
 \frac{1}{\epsilon} \int_{(t-\epsilon)^+}^t  \left( X_{s+\epsilon} -X_s \right) \indi_{\{ 0<X_{s+\epsilon} \}}ds 
 -
   \frac{1}{\epsilon} \int_{(t-\epsilon)^+}^t  \left( X_t -X_s \right) \indi_{\{ 0<X_{t} \}}ds . 
\end{eqnarray*}
Lemma \ref{articlerestezero} gives
\begin{equation}
\label{article0617}
 \sup_{t\in [0,  T]} | \Delta_2(t,  \epsilon) | \leqslant 2  C_\delta \epsilon^\delta. 
 \end{equation}
The relation (\ref{article0616}) leads us to study the convergence of $\widetilde{I}^2_{\epsilon}(t)$ to  $ X_t^+ +\frac{1}{2}L_t^0(X) $, as $\epsilon \to 0$. The idea of our approach is to express $\widetilde{I}^2_{\epsilon}(t)$   through $I^{4,1}_{\epsilon}(t) $  and $  I^{4,2}_{\epsilon}(t) $. First, we decompose $\widetilde{I}^2_{\epsilon}(t)$ as a sum of two terms:
$$
\widetilde{I}^2_{\epsilon}(t) 
 =  \frac{1}{\epsilon} \int_0^{t} X_{(s+\epsilon)\wedge t} \indi_{\{ 0<X_{(s+\epsilon) \wedge t }\}} ds - \frac{1}{\epsilon} \int_0^{t} X_{s} \indi_{\{ 0<X_{(s+\epsilon) \wedge t}\}} ds. 
$$
The change of variable $u=s+\epsilon$ in the first term gives:
\begin{eqnarray*}
\int_0^{t} X_{(s+\epsilon)\wedge t} \indi_{\{ 0<X_{(s+\epsilon) \wedge t }\}} ds &=&  \int_{\epsilon}^{t + \epsilon} X_{u\wedge t } \indi_{\{ 0<X_{u \wedge t}\}} du ,  \\
&=& \epsilon X_{ t } \indi_{\{ 0<X_{t}\}} +  \int_{0}^{t } X_{u } \indi_{\{ 0<X_{u }\}} du \\
&& - \int_{0}^{ \epsilon } X_{u } \indi_{\{ 0<X_{u}\}} du . 
\end{eqnarray*}
After easy calculations, we get:  
$$
\widetilde{I}^2_{\epsilon}(t)=X_t^+ 
+ \frac{1}{\epsilon} \int_{0}^{t} X_{u} \left( \indi_{\{ 0<X_u \}} - \indi_{\{ 0<X_{(u+\epsilon) \wedge t}\}} \right) du
-  \frac{1}{\epsilon} \int_{0}^{ \epsilon } X_u^+ du. 
$$
Next, we use again (\ref{article07cbis}), we obtain:
$$
\widetilde{I}^2_{\epsilon}(t)=X_t^+ 
+ I^{4,1}_{\epsilon}(t) +  I^{4,2}_{\epsilon}(t) -  \frac{1}{\epsilon} \int_{0}^{ \epsilon } X_u^+ du +  r^4_{\epsilon}(t),
$$ 
where it is recalled that $ I^{4,1}_{\epsilon}(t),   I^{4,2}_{\epsilon}(t), r^4_{\epsilon}(t)$ are defined by (\ref{article09b})-(\ref{article09breste}).  Hence,   we get
\begin{eqnarray}
\widetilde{I}^2_{\epsilon}(t)-  X_t^+ -\frac{1}{2}L_t^0(X) &=& \left( I^{4,1}_{\epsilon}(t) -\frac{1}{4}L_t^0(X) \right)  + \left( I^{4,2}_{\epsilon}(t) -\frac{1}{4}L_t^0(X) \right)\nonumber \\
&&-  \frac{1}{\epsilon} \int_{0}^{ \epsilon } X_u^+ du +  r^4_{\epsilon}(t).  \label{article0617b}
\end{eqnarray}
According to  Lemma \ref{article118}, we have:
\begin{equation}
\label{article0617c}
\sup_{t\in [0,  T]} \left|  \frac{1}{\epsilon} \int_{0}^{ \epsilon } X_u^+ du \right| \leqslant   \frac{C_\delta}{\delta + 1} \epsilon^\delta.
\end{equation}
Then, using  (\ref{article122}) and Theorem \ref{article11} comes to  
\begin{equation}
\label{article0618}
 \left\| \sup_{t\in [0,  T]} \left| \widetilde{I}^2_{\epsilon}(t)-  X_t^+ -\frac{1}{2}L_t^0(X)\right| \right\|_{L^2(\Omega)}  \leqslant  C \epsilon^{\frac{\delta}{2}}.  
\end{equation}
It is clear that Theorem \ref{article06b} is a direct consequence of  (\ref{article122}),(\ref{article0613}), (\ref{article068brownien}), (\ref{article0619}), (\ref{article0615}), (\ref{article0617}), 
(\ref{article0617c}) and (\ref{article0618}).   \qed

   \textbf{4.}  We now demonstrate Proposition \ref{article06bb}. Let us consider a sequence $ (\epsilon_n)_{n\in \N}$  of positive real numbers decreasing to 0,   such that $\sum_{i=1 }^{\infty} \sqrt{\epsilon_i} <\infty$.

  First, we show  the almost sure convergence of $I^{4,2}_{\epsilon_n} (t)$ to $ \frac{1}{4}L^0_t(X)$, as $n\to \infty$. Note that Theorems \ref{article11} and \ref{article06b} do not permit to obtain the a.s. convergence results stated in Proposition  \ref{article06bb}, via the Borel Cantelli lemma, since we cannot take $\delta=\frac{1}{2}$. \\
Recall the identity (\ref{article123}):
$$I^{4,2}_\epsilon(t) -   \frac{1}{4}L^0_t(X) = D^1_{\epsilon}(t) + D^2_{\epsilon}(t)+ D^3_{\epsilon}(t) +  \left(  \frac{1}{\epsilon} \int_{0}^t X_{u}^+ \Phi \left( -\frac{X_u}{\sqrt{ \epsilon}}\right) - \frac{1}{4}L^0_t(X) \right).$$
According to (\ref{article121}), (\ref{article115}) and (\ref{article116}) , the quantities \\
$ \left(  \frac{1}{\epsilon} \int_{0}^t X_{u}^+ \Phi \left( -\frac{X_u}{\sqrt{ \epsilon}}\right) - \frac{1}{4}L^0_t(X) \right)$, $D^2_{\epsilon}(t)$ and $D^3_{\epsilon}(t)$ goes to 0, as $\epsilon \to 0$, with a rate of decay of order $\delta < \frac{1}{2}$. But it does not matter since the convergence holds in the almost sure sense.

   We now focus on $D^1_{\epsilon}(t) $. According to (\ref{article120}), we have
$E\left[ \sup_{t\in [0,  T]} \left( D^1_{\epsilon}(t)  \right)^2\right] \leqslant C \sqrt{\epsilon}$.  Then,   the Borel Cantelli Lemma implies that,   for all $t\geqslant 0,   \left( D^1_{\epsilon_n}(t)\right)_{n\in\N}$ converge almost surely  uniformly on $[0,  T]$.   Hence,   we get the a.s.  convergence of $ I^{4,2}_{\epsilon_n}(t)$ to $\frac{1}{4}L^0_t(X)$, as $n\to\infty$.

The convergence of $I^{4,1}_{\epsilon_n} (t)$ to $\frac{1}{4}L^0_t(X)$ may be obtained by using the symmetry of Brownian motion.

  Let us now investigate the convergence of  $J_{\epsilon_n} (t)$. It is clear that (\ref{article0613}), (\ref{article068}), (\ref{article0616}) and (\ref{article0617b}) imply that $J_{\epsilon}(t) - L^0_t(X)$ is equal to:
$$
\begin{array}{c}
\widetilde{I}^1_\epsilon(t) + \Delta^1_\epsilon(t) + \left(I^{4,1}_\epsilon(t) -   \frac{1}{4}L^0_t(X) \right)+  \left(I^{4,2}_\epsilon(t) -   \frac{1}{4}L^0_t(X) \right)   \\
 + \Delta^2_\epsilon(t) +r^4_\epsilon (t) -\frac{1}{\epsilon}\int_0^\epsilon X_u^+ du.
\end{array}
$$
%
We proceed as in the convergence of $ I^{4,2}_{\epsilon_n}(t)$. From inequalities
(\ref{article0619}), (\ref{article0617}), (\ref{article122}) and (\ref{article0617c}), it may be deduced that  $ \Delta^1_\epsilon(t) + \Delta^2_\epsilon(t) + r^4_\epsilon (t) -\frac{1}{\epsilon}\int_0^\epsilon X_u^+ du $ tends a.s to 0 as $\epsilon \to 0$. Note that we have already shown that  $ \left(I^{4,i}_{\epsilon_n}(t) -   \frac{1}{4}L^0_t(X) \right)$, $i=1,2$, converges a.s. as $n\to \infty$. Finally, the a.s. convergence of $\widetilde{I}^1_{\epsilon_n}(t)$, as $n\to \infty$, may be obtained through  (\ref{article0615}).
\qed

{\small
\def\refname{References}
\bibliography{bibthese}
\bibliographystyle{plain}
\nocite{*}
}

\end{document}